\documentclass[onefignum,onetabnum]{siamart190516}

\usepackage{lipsum}
\usepackage{amsfonts}
\usepackage{graphicx}
\usepackage{epstopdf}
\usepackage{algorithmic}
\usepackage{geometry}
\usepackage{amssymb}
\usepackage{caption}
\usepackage{subcaption}
\usepackage{enumitem}
\usepackage{mathtools}
\usepackage[all]{xy}
\usepackage{xcolor}
\usepackage{soul}
\usepackage{cite}

\graphicspath{{figures/}}

\ifpdf
\DeclareGraphicsExtensions{.eps,.pdf,.png,.jpg}
\else
\DeclareGraphicsExtensions{.eps}
\fi

\newsiamremark{remark}{Remark}
\newsiamremark{hypothesis}{Hypothesis}
\crefname{hypothesis}{Hypothesis}{Hypotheses}
\newsiamthm{claim}{Claim}

\headers{Adaptive Variational Integrators and Symplectic Accelerated Optimization}{Valentin Duruisseaux, Jeremy Schmitt, and Melvin Leok}

\title{Adaptive Hamiltonian Variational Integrators and Applications to Symplectic Accelerated Optimization}

\author{Valentin Duruisseaux, Jeremy Schmitt, and Melvin Leok}

\usepackage{amsopn}

\begin{document}

\maketitle

% REQUIRED
\begin{abstract}
It is well known that symplectic integrators lose their near energy preservation properties when variable time-steps are used. The most common approach to combine adaptive time-steps and symplectic integrators involves the Poincar\'e transformation of the original Hamiltonian. In this article, we provide a framework for the construction of variational integrators using the Poincar\'e transformation. Since the transformed Hamiltonian is typically degenerate, the use of Hamiltonian variational integrators based on Type II or Type III generating functions is required instead of the more traditional Lagrangian variational integrators based on Type I generating functions. Error analysis is provided and numerical tests based on the Taylor variational integrator approach of \cite{ScShLe2017} to time-adaptive variational integration of Kepler's 2-Body problem are presented. Finally, we use our adaptive framework together with the variational approach to accelerated optimization presented in \cite{WiWiJo16} to design efficient variational and non-variational explicit integrators for symplectic accelerated optimization. 
\end{abstract}

% REQUIRED

% REQUIRED

\section{Introduction}

Symplectic integrators form a class of geometric numerical integrators of interest since, when applied to Hamiltonian systems, they yield discrete approximations of the flow that preserve the symplectic 2-form (see \cite{HaLuWa2006}). The preservation of symplecticity results in the preservation of many qualitative aspects of the underlying dynamical system. In particular, when applied to conservative Hamiltonian systems, symplectic integrators show excellent long-time near-energy preservation. However, when symplectic integrators were first used in combination with variable time-steps, the near-energy preservation was lost and the integrators performed poorly (see \cite{CalSan93, GlaDunCan91}). Backward error analysis provided justification both for the excellent long-time near-energy preservation of symplectic integrators and for the poor performance experienced when using variable time-steps (see Chapter IX of \cite{HaLuWa2006}). Backward error analysis shows that symplectic integrators can be associated with a modified Hamiltonian in the form of a powers series in terms of the time-step. The use of a variable time-step results in a different modified Hamiltonian at every iteration where the time-step is changed, which is the source of the poor energy conservation. There has been a great effort to circumvent this problem, and there have been many successes. However, there has yet to be a unified general framework for constructing adaptive symplectic integrators. In this paper, we contribute to this effort by demonstrating how Hamiltonian variational integrators \cite{LeZh2011} can be used to systematically construct symplectic integrators that allow for the use of variable time-steps.

The use of variable time-steps is motivated by the observation that the global error estimates for a numerical method depend in part on the maximum local truncation error, and this in turn is related to both the time-step and the magnitude of the $(r+1)$-derivatives of the solution for a $r$-order numerical method. For a fixed number of time-steps, the maximum local truncation error is minimized if the local truncation error is equidistributed over the time intervals. In turn, this can be achieved if, for example, the time-step is chosen to be an appropriate function of the reciprocal of the relevant derivative of the solution. This derivative can be estimated \textit{a posteriori} by comparing methods with different orders of accuracy, or methods with the same order of accuracy but different error constants. Alternatively, in the Kepler two-body problem, for example, Kepler's second law states that the line joining the planet and the Sun sweeps out equal areas during equal intervals of time, so the angular velocity of the planet is proportional to the reciprocal of the radius squared, which gives an \textit{a priori} bound. In essence, variable time-steps are chosen to control the error incurred at each time-step, which in turn affects the global accuracy of the numerical trajectory.

The goal of this paper is to develop an analogue of the methods derived using the framework of \cite{Ha1997, Re1999}, but directly in terms of generating functions of symplectic maps. These prior results are based on symplectic (partitioned) Runge--Kutta methods, which are related to Type I generating functions~\cite{Su1990}, but we desire an explicit characterization of the flow maps of time-adaptive Hamiltonian systems so that we can employ the Hamiltonian variational integrator framework instead.

Variational integrators provide a systematic method for constructing symplectic integrators of arbitrarily high-order based on the discretization of Hamilton's principle~\cite{MaWe2001, HaLe2012}, or equivalently, by the approximation of generating functions, but there has not been a systematic attempt to incorporate time-adaptivity into the setting of variational integrators. This is due to the fact that the Poincar\'e transformed Hamiltonian that is used is in general degenerate, so there is no corresponding Lagrangian analogue, which prevents the use of traditional variational integrators that are based on a Lagrangian formulation of mechanics and involve the construction of a discrete Lagrangian that approximates a Type I generating function given by Jacobi's solution of the Hamilton--Jacobi equation. Instead, we propose the use of Hamiltonian variational integrators~\cite{LeZh2011}, which are based on Type II and Type III generating functions that have no difficulty with this degeneracy.

After a brief introduction to variational integrators in Section \ref{SectionBasicTheory}, we will review the construction of Type II and Type III Hamiltonian Taylor variational integrators from \cite{ScShLe2017}  and present a new theorem concerning their order of accuracy in Section \ref{HTVISection}. We will then present a framework for variable time-step variational integrators in Section \ref{PoincareSubsection}, derive corresponding error analysis results in Section \ref{ErrorAnalysisSection}, and test our approach with Hamiltonian Taylor variational integrators on Kepler's 2-Body problem in Section \ref{KeplerSection}. Finally, in Section \ref{OptimizationSection}, we will design efficient variational and non-variational explicit integrators for symplectic accelerated optimization, using our adaptive approach in the variational framework to accelerated optimization introduced in \cite{WiWiJo16}.  \\

\section{Hamiltonian Variational Integrators}\label{SectionHamiltonianVI}

\subsection{Variational Integration} \label{SectionBasicTheory}

Variational integrators are derived by discretizing Hamilton's principle, instead of discretizing Hamilton's equations directly. As a result, variational integrators are symplectic, preserve many invariants and momentum maps, and have excellent long-time near-energy preservation (see \cite{MaWe2001}). \\

\noindent\underline{\textbf{Type I:}} Traditionally, variational integrators have been designed based on the Type I generating function known as the discrete Lagrangian, $L_d:Q \times Q \mapsto \mathbb{R}$. The exact discrete Lagrangian of the true flow of Hamilton's equations can be represented in both a variational form and in a boundary-value form. The latter is given by
\begin{equation}
	L_d^E(q_0,q_1;h)=\int_0^h L(q(t),\dot q(t)) dt  \label{exact_Ld},
\end{equation}
where $q(0)=q_0,$ $q(h)=q_1,$ and $q$ satisfies the Euler--Lagrange equations over the time interval $[0,h]$. A variational integrator is defined by constructing an approximation $L_d:Q \times Q \mapsto \mathbb{R}$ to $L_d^E$, and then applying the discrete Euler--Lagrange equations,
\begin{equation}
	p_k=-D_1 L_d(q_k, q_{k+1}),\qquad p_{k+1}=D_2 L_d(q_k, q_{k+1}),  \label{IDEL}
\end{equation}
which implicitly define the integrator $\tilde{F}_{L_d}:(q_k,p_k)\mapsto(q_{k+1},p_{k+1})$, where $D_i$ denotes a partial derivative with respect to the $i$-th argument. The error analysis is greatly simplified via Theorem 2.3.1 of \cite{MaWe2001}, which states that if a discrete Lagrangian, $L_d:Q\times Q\rightarrow\mathbb{R}$, approximates the exact discrete Lagrangian $L_d^E:Q\times Q\rightarrow\mathbb{R}$ to order $r$, i.e.,
\begin{equation} 
	L_d(q_0, q_1;h)=L_d^E(q_0,q_1;h)+\mathcal{O}(h^{r+1}) ,
\end{equation}
then the discrete Hamiltonian map $\tilde{F}_{L_d}:(q_k,p_k)\mapsto(q_{k+1},p_{k+1})$, viewed as a one-step method, has order of accuracy $r$. Many other properties of the integrator, such as momentum conservation properties of the method, can be determined by analyzing the associated discrete Lagrangian, as opposed to analyzing the integrator directly. 

More recently, variational integrators have been extended to the framework of Type II and Type III generating functions, commonly referred to as discrete Hamiltonians (see \cite{LaWe2006, LeZh2011,ScLe2017}). Hamiltonian variational integrators are derived by discretizing Hamilton's phase space principle. \\

\noindent \underline{\textbf{Type II:}} The boundary-value formulation of the exact Type II generating function of the time-$h$ flow of Hamilton's equations is given by the exact discrete right Hamiltonian,
\begin{equation}
	H_d^{+,E}(q_0,p_1;h) =  p_1^\top q_1 - \int_0^h \left[ p(t)^\top \dot{q}(t)-H(q(t), p(t)) \right] dt, \label{exact_Hd}
\end{equation}
where $(q,p)$ satisfies Hamilton's equations with boundary conditions $q(0)=q_0$ and $p(h)=p_1$. A Type II Hamiltonian variational integrator is constructed by using an approximate discrete Hamiltonian $H_d^+$, and applying the discrete right Hamilton's equations,
\begin{equation}\label{Discrete Right Eq}
	p_0=D_1H_d^+(q_0,p_1), \qquad q_1=D_2H_d^+(q_0,p_1),
\end{equation}
which implicitly defines the integrator, $\tilde{F}_{H_d^+}:(q_0,p_0) \mapsto (q_1,p_1)$.

Theorem 2.3.1 of \cite{MaWe2001}, which simplified the error analysis for Lagrangian variational integrators, has an analogue for Hamiltonian variational integrators. Theorem 2.2 in \cite{ScLe2017} states that if a discrete right Hamiltonian $H^+_d$ approximates the exact discrete right Hamiltonian $H_d^{+,E}$ to order $r$, i.e.,
\begin{equation} 
	H^+_d(q_0, p_1;h)=H_d^{+,E}(q_0,p_1;h)+\mathcal{O}(h^{r+1}),
\end{equation}
then the discrete right Hamilton's map $\tilde{F}_{H^+_d}:(q_k,p_k)\mapsto(q_{k+1},p_{k+1})$, viewed as a one-step method, is order $r$ accurate. \\

\noindent \underline{\textbf{Type III:}}  The boundary-value formulation of the exact Type III generating function of the time-$h$ flow of Hamilton's equations is given by the exact discrete left Hamiltonian,
\begin{equation}
	H_d^{-,E}(q_1,p_0;h) =  - p_0^\top q_0  - \int_0^h \left[ p(t)^\top \dot{q}(t)-H(q(t), p(t)) \right] dt \label{exact_LeftHd},
\end{equation}
where $(q,p)$ satisfies Hamilton's equations with boundary conditions $q(h)=q_1$ and $p(0)=p_0$. A Type III Hamiltonian variational integrator is constructed by using an approximate discrete left Hamiltonian $H_d^-$, and applying the discrete left Hamilton's equations,
\begin{equation} \label{Discrete Left Eq}
	p_1= - D_1H_d^-(q_1,p_0), \qquad q_0 = - D_2H_d^-(q_1,p_0),
\end{equation}
which implicitly defines the integrator, $\tilde{F}_{H_d^-}:(q_0,p_0) \mapsto (q_1,p_1)$. As mentioned in \cite{ScLe2017}, the proof of Theorem 2.2 in \cite{ScLe2017} can be easily adjusted to prove an equivalent theorem for the discrete left Hamiltonian case, which states that if a discrete left Hamiltonian $H^-_d$ approximates the exact discrete left Hamiltonian $H_d^{-,E}$ to order $r$, i.e.,
\begin{equation} 
	H^-_d(q_1, p_0;h)=H_d^{-,E}(q_1,p_0;h)+\mathcal{O}(h^{r+1}),
\end{equation}
then the discrete left Hamilton's map $\tilde{F}_{H^-_d}:(q_k,p_k)\mapsto(q_{k+1},p_{k+1})$, viewed as a one-step method, is order $r$ accurate.  \\

Examples of Hamiltonian variational integrators include Galerkin variational integrators ~\cite{LeZh2011}, Prolongation-Collocation variational integrators~\cite{LeSh2011}, and Taylor variational integrators~\cite{ScShLe2017}. %We present new results about the order of accuracy of Hamiltonian Taylor variational integrators in the next section.
In many cases, the Type I and Type II/III approaches will produce equivalent integrators. This equivalence has been established in \cite{ScShLe2017} for Taylor variational integrators provided the Lagrangian is hyperregular, and in \cite{LeZh2011} for generalized Galerkin variational integrators constructed using the same choices of basis functions and numerical quadrature formula provided the  Hamiltonian is hyperregular. However, Hamiltonian and Lagrangian variational integrators are not always equivalent. In particular, it was shown in \cite{ScLe2017} that even when the Hamiltonian and Lagrangian integrators are analytically equivalent, they might still have different numerical properties because of numerical conditioning issues. Even more to the point, Lagrangian variational integrators cannot always be constructed when the underlying Hamiltonian is degenerate, and in that situation, Hamiltonian variational integrators are the more natural choice. Depending on the form of the Hamiltonian and the method used to design the corresponding approximate discrete Hamiltonian, one of the Type II or Type III approaches might be more convenient than the other, in the sense that it might allow for an explicit algorithm or might allow for higher-order methods given some constraints on the type of methods permitted. In Section \ref{AdaptiveSection}, we will examine a transformation commonly used to construct variable time-step symplectic integrators, which results in a degenerate Hamiltonian in most cases of interest, such as the optimization application considered in Section \ref{OptimizationSection}. We will apply Hamiltonian variational integrators to the resulting transformed Hamiltonian system. For the optimization application presented in Section \ref{OptimizationSection}, we will prefer Type II Hamiltonian Taylor variational integrators to their Type III analogues, and this choice will be justified carefully based on the order and explicitness of the resulting methods. \\

\subsection{Hamiltonian Taylor Variational Integrators (HTVIs)} \label{HTVISection}

We now present Hamiltonian Taylor variational integrators \cite{ScShLe2017}, together with a new theorem concerning their order of accuracy, which is analogous to Theorem 3.1 in \cite{ScShLe2017} for their Lagrangian counterpart. A discrete approximate Hamiltonian is constructed by approximating the flow map and the trajectory associated with the boundary values using a Taylor method, and approximating the integral by a quadrature rule. The HTVI is then generated by the discrete Hamilton's equations associated to that discrete Hamiltonian. More explicitly, we first construct the $r$-order and $(r+1)$-order Taylor methods $\Psi_h^{(r)}$ and $\Psi_h^{(r+1)}$ approximating the exact time-$h$ flow map $\Phi _h : T^*Q \rightarrow T^*Q$ corresponding to Hamilton's equation $ \dot{z}= \varphi (z)$,  where $z = (q,p)$:

\begin{equation}
	\Psi_h^{(r)}(z_0) = z_0 + \sum_{k=1}^{r}{\frac{h^k}{k!} \varphi^{(k-1)}(z_0})  .
\end{equation}

Let $\pi_{T^*Q}:(q,p)\mapsto p$ and $\pi_{Q}:(q,p)\mapsto q$. Given a quadrature rule of order $s$ with weights and nodes $(b_i,c_i)$ for $i=1,...,m$, the Type II and Type III integrators are then constructed as follows:

\begin{table}[!h]
	\resizebox{15cm}{!} {
	\begin{tabular}{l|l}
		\textbf{Type II HTVI} & \textbf{Type III HTVI}  \\ \hline
		Approximate $p(0)$ by the solution $\tilde{p}_0$ of  & Approximate $q(0)$ by the solution $\tilde{q}_0$ of  \\  
		$ \qquad \qquad  p_1 = \pi_{T^*Q} \circ \Psi_h^{(r)}(q_0,\tilde{p}_0) .$ & $ \qquad \qquad q_1 = \pi_{Q} \circ \Psi_h^{(r+1)}(\tilde{q}_0,p_0). $ \\
		Approximate $(q(c_i h),p(c_i h))$ via & Approximate $(q(c_i h),p(c_i h))$ via \\  $ \qquad \qquad (q_{c_i},p_{c_i})  = \Psi_{c_i h}^{(r)}(q_0,\tilde{p}_0) .$ &  $ \qquad \qquad (q_{c_i},p_{c_i})  = \Psi_{c_i h}^{(r)}(\tilde{q}_0,p_0).$  \\ Approximate $q_1$ via & \\ $ \qquad \qquad \tilde{q}_1 = \pi_{Q}  \circ \Psi_h^{(r+1)}(q_0,\tilde{p}_0).$ \\  Continuous Legendre transform:  $\dot{q}_{c_i} = \frac{\partial H}{\partial p_{c_i}}$. &  Continuous Legendre transform:  $\dot{q}_{c_i} = \frac{\partial H}{\partial p_{c_i}}$. \\ 
		Then, the quadrature rule gives & Then, the quadrature rule gives  \\
		\small $H_d^+(q_0,p_1)  = p_1^\top \tilde{q}_1 -  h \sum_{i=1}^{m}{b_i \left[  p_{c_i}^\top \dot{q}_{c_i} - H(q_{c_i},p_{c_i})  \right]} $ & \small $ H_d^-(q_1,p_0) = -p_0^\top \tilde{q}_0 - h \sum_{i=1}^{m}{b_i \left[  p_{c_i}^\top \dot{q}_{c_i} - H(q_{c_i},p_{c_i})  \right]} $ \\
		The discrete right Hamilton's equations \eqref{Discrete Right Eq}  & The discrete left Hamilton's equations  \eqref{Discrete Left Eq}  \\ then define the variational integrator. & then define the variational integrator.
	\end{tabular}
}
\end{table}

Taylor variational integrators were inspired by a resurgence of interest in high-order Taylor methods for celestial mechanics that has been fueled by the continued progress in Automatic Differentiation software (see \cite{em/1120145574,10.1093/mnras/stab1032,BARRIO2005525,Pearlmutter07lazymultivariate,Bettencourt2019TaylorModeAD,Neidinger2005DfC,Neidinger2013Err}). Implicit modified Taylor methods have been proposed to deal with stiff ODEs \cite{osti_6097907}, while Taylor variational integrators provide a class of Taylor-based integrators to deal with conservative Hamiltonian systems, and can be viewed as a predictor-corrector method that applies a symplectic correction to the Taylor method. For high-order Taylor methods, the key to an efficient implementation relies upon efficient Automatic Differentiation software to compute higher-order gradients. 

We now present a theorem specifying the order of accuracy of the resulting Hamiltonian Taylor variational integrators:

\begin{theorem} \label{HTVITheorem}
	If the Hamiltonian $H$ and $\frac{\partial H}{\partial p}$ are Lipschitz continuous in both variables, then the discrete Hamiltonian $H_d^{\pm}$ obtained using the above construction, approximates $H_d^{\pm,E}$ with at least order of accuracy $\min{(r+1,s)}$. By Theorem 2.2 in \cite{ScLe2017} (or its analogue for the left Hamiltonian case), the associated discrete Hamiltonian map has the same order of accuracy.
	\begin{proof}
		See Appendix \ref{HTVIProof}.  
	\end{proof} 
\end{theorem} 

\hfill 

\section{Adaptive Integrators and Variational Error Analysis} \label{AdaptiveSection}

\subsection{The Poincar\'e Transformation and Discrete Hamiltonians} \label{PoincareSubsection}

Given an autonomous Hamiltonian $H(q,p)$, and a desired transformation of time $t \mapsto \tau$ described by the monitor function $g(q,p)$ via
\begin{equation}
	\frac{dt}{d\tau} = g(q,p),
\end{equation}
a new Hamiltonian system is constructed using the Poincar\'e transformation,
\begin{equation} \label{TransformedH}
	\bar{H}(\bar{q},\bar{p}) = g(q,p) \left(H(q,p) + p^t \right),
\end{equation}
where $\bar{q} = \begin{bmatrix} q \\ q^t \end{bmatrix} $ and $\bar{p} = \begin{bmatrix} p \\ p^t \end{bmatrix} $. We will make the common choice $q^t=t$ and $p^t=-H(q(0),p(0))$, so that $\bar{H}(\bar{q},\bar{p})=0$ along all integral curves through $(\bar{q}(0),\bar{p}(0))$. The time $t$ shall be referred to as the physical time, while $\tau$ will be referred to as the fictive time. 

In general, along an integral curve through $(\bar{q}(0),\bar{p}(0))$, 
\begin{equation} \label{d2Hdp2}
	\frac{\partial^2 \bar{H}}{\partial \bar{p}^2} = \begin{bmatrix} \frac{\partial H}{\partial p}\nabla_pg(q,p)^\top+g(q,p)\frac{\partial^2 H}{\partial p^2}+\nabla_pg(q,p)\frac{\partial H}{\partial p}^\top  \ \ \nabla_pg(q,p) \\ \nabla_pg(q,p)^\top \qquad \qquad \qquad \qquad \qquad \qquad \qquad \qquad 0
	\end{bmatrix},
\end{equation}
which can be singular for many initial Hamiltonians $H$ and choices of monitor function $g$.

Most of the prior literature on variable time-step symplectic integrators cited in this paper focus exclusively on monitor functions that are only a function of position, in which case $\frac{\partial^2 \bar{H}}{\partial \bar{p}^2}$ is singular, and the associated Legendre transformation, $\mathbb{F}\bar{H}:T^*Q \rightarrow TQ$ is non-invertible, which is to say that the resulting transformed Hamiltonian is degenerate and there is no corresponding Lagrangian formulation. Therefore, the Type II and Type III Hamiltonian variational integrator frameworks are the most general and natural way to derive variable time-step variational integrators. 

The exact Type II generating function for the transformed Hamiltonian is given by
\begin{equation}
	\bar{H}^{+,E}_d(\bar{q}_0,\bar{p}_1;h) = \bar{p}_1^\top \bar{q}_1 - \int_0^h \left( \bar{p}(\tau)^\top \dot{\bar{q}}(\tau) - \bar{H}(\bar{q}(\tau),\bar{p}(\tau)) \right) d\tau, \label{trans_exact_Hd}
\end{equation}
where $(\bar{q}(\tau), \bar{p}(\tau))$ satisfy the Hamilton's equations corresponding to the Poincar\'e transformed Hamiltonian $\bar{H}$, with boundary conditions $\bar{q}(0) = \bar{q}_0$ and  $\bar{p}(h) = \bar{p}_1$. This exact discrete right Hamiltonian implicitly defines a symplectic map with respect to the symplectic form $\bar{\omega}(\bar{p}_k,\bar{q}_k)$ on $T^*\bar{Q}$ via the discrete Legendre transforms given by
\begin{equation}
	\bar{p}_0 = \frac{\partial \bar{H}_d^{+,E}}{\partial \bar{q}_0}, \qquad \bar{q}_1 = \frac{\partial \bar{H}_d^{+,E}}{\partial \bar{p}_1}.
\end{equation}

Similarly, the exact Type III generating function for the transformed Hamiltonian is given by
\begin{equation}
	\bar{H}^{-,E}_d(\bar{q}_0,\bar{p}_1;h) = -\bar{p}_0^\top\bar{q}_0 - \int_0^h \left( \bar{p}(\tau)^\top\dot{\bar{q}}(\tau) - \bar{H}(\bar{q}(\tau),\bar{p}(\tau)) \right) d\tau, \label{lefttrans_exact_Hd}
\end{equation}
where $(\bar{q}(\tau), \bar{p}(\tau))$ satisfy the Hamilton's equations corresponding to the Poincar\'e transformed Hamiltonian $\bar{H}$ with boundary conditions $\bar{q}(h) = \bar{q}_1$ and $\bar{p}(0) = \bar{p}_0$. This exact discrete left Hamiltonian implicitly defines a symplectic map with respect to the symplectic form $\bar{\omega}(\bar{p}_k,\bar{q}_k)$ on $T^*\bar{Q}$ via the discrete Legendre transforms given by
\begin{equation}
	\bar{p}_1 = - \frac{\partial \bar{H}_d^{-,E}}{\partial \bar{q}_1}, \qquad \bar{q}_0 = - \frac{\partial \bar{H}_d^{-,E}}{\partial \bar{p}_0}.
\end{equation}

Our approach is to construct Hamiltonian variational integrators using a discrete Hamiltonian $\bar{H}^{\pm}_d$ that approximates the corresponding exact discrete Hamiltonian $\bar{H}^{\pm,E}_d$ to order $r$. The resulting integrator will be symplectic with constant time-step in fictive time $\tau$ and more importantly with the desired variable time-step in physical time $t$ via $\frac{dt}{d\tau} = g(q,t,p)$. It is important to note that this method will be symplectic in two different ways. It will be symplectic both with respect to the symplectic form $d\bar{p} \wedge d\bar{q}$ and with respect to the symplectic form $dp \wedge dq$. Since $\dot{p}^t =0$, $p^t$ is constant and $dp^t_{k} \wedge dq^t_{k} =0$, the symplectic form in generalized coordinates is given by
\begin{equation}
	\bar{\omega}(\bar{p}_k,\bar{q}_k) = d\bar{p}_k \wedge d\bar{q}_k 
	= \sum_{i=1}^{n+1}  d\bar{p}_{k,i} \wedge d\bar{q}_{k,i} 
	= \sum_{i=1}^{n}  dp_{k,i} \wedge dq_{k,i} = \omega(p_k,q_k).
\end{equation}

A symplectic variable time-step method was proposed independently in \cite{Ha1997} and \cite{Re1999}, which applied a symplectic integrator to the Poincar\'e transformed Hamiltonian. In \cite{Ha1997}, it is noted that one of the first applications of the Poincar\'e transformation was by Levi-Civita, who applied it to the three-body problem. A more in-depth discussion of such time transformations can be found in \cite{Struc2005}. Further work using this type of transformation has been published, such as \cite{BlBu2004, BlIs2012} which focused on developing symplectic, explicit, splitting methods with variable time-steps.  

The novelty of our approach consists in discretizing the Type II or Type III generating function for the flow of Hamilton's equations, where the Hamiltonian is given by the Poincar\'e transformation. Therefore, we are constructing variational integrators, and in particular Hamiltonian variational integrators (see \cite{LaWe2006, LeZh2011}). The use of Type II or Type III integrators is justified by the degeneracy of the Hamiltonian which implies that there is no corresponding Type I Lagrangian formulation. This approach works seamlessly with existing methods and theorems of Hamiltonian variational integrators, but now the system under consideration is the transformed Hamiltonian system resulting from the Poincar\'e transformation. The methods of \cite{Ha1997, Re1999} include the possibility of applying a given variational integrator to the transformed differential equations. Our approach gives a framework for constructing variational integrators at the level of the generating function by using the Poincar\'e transformed discrete right Hamiltonian.

\begin{remark*}
	Other approaches to variable time-step variational integrators can be found in \cite{KaMaOr99, FuMo2006, Nair2012}. In particular, \cite{KaMaOr99} is inspired by the result of \cite{GeMa88}, which states that constant time-step symplectic integrators of autonomous Hamiltonian systems cannot exactly conserve the energy unless it agrees with the exact flow map up to a time reparametrization. Therefore, they sought a variable time-step energy-conserving symplectic integrator in an expanded non-autonomous system. However, symplecticity is with respect to the space-time symplectic form $dp \wedge dq + dH \wedge dt$. The time-step is determined by enforcing discrete energy conservation, which arises as a consequence of the fact that energy is the Noether quantity associated with time translational symmetry. An extended Hamiltonian is used, similar in spirit to the Poincar\'e transformation. An approach that builds off this idea and space-time symplecticity was presented in \cite{Nair2012}, and a less constrained choice of time-step was allowed. In \cite{FuMo2006}, adaptive variational integrators are constructed using a transformation of the Lagrangian, which is motivated by the Poincar\'e transformation, but it is not equivalent. The lack of equivalence is not surprising, since the Poincar\'e transformed Hamiltonian is degenerate for their choice of monitor functions. As a consequence, the phase space path is not preserved.
\end{remark*}

Note that our framework can be extended to the case where the original Hamiltonian $H$ and the chosen monitor function $g$ depend explicitly on time $t$ (inspired by \cite{Ha1997}). Given a time-dependent Hamiltonian $H(q,t,p)$, consider a desired  transformation of time $t \mapsto \tau$, given by $\frac{dt}{d\tau} = g(q,t,p).$ Then, we can define $\bar{q} = \begin{bmatrix} q \\ q^t \end{bmatrix} $ where $q^t = t$ and $\bar{p} = \begin{bmatrix} p \\ p^t \end{bmatrix} $ where $p^t$ is the conjugate momentum for $q^t = t$ with $p^t(0) = - H(q(0),0,p(0))$. Consider the new Hamiltonian system given by the Poincar\'e transformation
\begin{equation}
	\bar{H}(\bar{q},\bar{p}) = g(q,q^t,p) \left(H(q,q^t,p) + p^t \right).
\end{equation}
The corresponding equations of motion in the extended phase space are then given by
\begin{equation}
	\dot{\bar{q}}  = \frac{\partial \bar{H}}{\partial \bar{p}},  \qquad\qquad  \dot{\bar{p}} =  -\frac{\partial \bar{H}}{\partial \bar{q}}.    
\end{equation}
Suppose $(\bar{Q}(\tau) , \bar{P}(\tau))$ are solutions to these extended equations of motion, and let $(q(t),p(t))$ solve Hamilton's equations for the original Hamiltonian $H$. Then
\begin{equation}
	\bar{H}(\bar{Q}(\tau) , \bar{P}(\tau))  = \bar{H}(\bar{Q}(0) , \bar{P}(0)) = 0.
\end{equation}
Thus, the components $(Q(\tau) , P(\tau))$ in the original phase space of the solutions $(\bar{Q}(\tau) , \bar{P}(\tau))$ satisfy
\begin{equation}H(Q(\tau) , \tau , P(\tau)) = - p^t(\tau), \qquad  H(Q(0) , 0, P(0)) = - p^t(0) = H(q(0),0,p(0)) . \end{equation}
Then,  $(Q(\tau) , P(\tau))$ and $(q(t),p(t))$ both satisfy Hamilton's equations for the original Hamiltonian $H$ with the same initial values, so they must be the same. As before,
\begin{equation}
	\frac{\partial^2 \bar{H}}{\partial \bar{p}^2} = \begin{bmatrix} \frac{\partial H}{\partial p}\nabla_pg(\bar{q},p)^\top +g(\bar{q},p)\frac{\partial^2 H}{\partial p^2}+\nabla_pg(\bar{q},p)\frac{\partial H}{\partial p}^\top  \ \ \nabla_pg(\bar{q},p) \\ \nabla_pg(\bar{q},p)^\top  \qquad \qquad \qquad \qquad \qquad \qquad \qquad \qquad 0
	\end{bmatrix},
\end{equation}
will be singular in many cases, so Hamiltonian variational integrators are the most general and natural way to derive variable time-step variational integrators.  \\

\subsection{Variational Error Analysis} \label{ErrorAnalysisSection}

The standard error analysis for Hamiltonian variational integrators assumes a non-degenerate Hamiltonian, i.e., $\text{det}\left(\frac{\partial^2 \bar{H}}{\partial \bar{p}^2}\right) \neq 0$ (see \cite{ScLe2017}), which might not be the case for the Poincar\'e transformed Hamiltonian. The non-degeneracy of the Hamiltonian ensures that we can apply the usual implicit function theorem to the discrete Hamilton's equations, and the proof of the standard error analysis theorem relies upon Lemma 2.3 of \cite{ScLe2017}:

\begin{lemma} \label{InvertibilityLemma}
	Let $f_1,g_1,e_1,f_2,g_2,e_2 \in C^r$ ($r$-times continuously differentiable) be such that
	\begin{equation}f_1(x,h)=g_1(x,h)+h^{r+1}e_1(x,h), \qquad f_2(x,h)=g_2(x,h)+h^{r+1}e_2(x,h).\end{equation}
	Then, there exists functions $e_{12}$ and $\bar{e}_1$ bounded on compact sets such that
	\begin{equation}f_2(f_1(x,h),h)=g_2(g_1(x,h),h)+h^{r+1}e_{12}(g_1(x,h),h),\end{equation}
	\begin{equation}f_1^{-1}(y)=g_1^{-1}(y)+h^{r+1}\bar{e}_1(y).\end{equation}
\end{lemma} 

Given a discrete Hamiltonian $H_d^\pm$, we introduce the discrete fiber derivatives (or discrete Legendre transforms), $\mathbb{F}^\pm H_d^\pm$,
\begin{equation*}
\begin{aligned}
	\mathbb{F}^+H_d^+(q_0,p_1)&: (q_0,p_1)\mapsto (D_2 H_d^+(q_0,p_1),p_1),  \quad  \mathbb{F}^+H_d^-(q_1,p_0)&: (q_1,p_0)\mapsto (q_1,-D_1 H_d^-(q_1,p_0)),\\
	\mathbb{F}^-H_d^+(q_0,p_1)&: (q_0,p_1)\mapsto (q_0 , D_1 H_d^+ (q_0,p_1)), \quad \mathbb{F}^-H_d^-(q_1,p_0)&: (q_1,p_0)\mapsto (-D_2 H_d^-(q_1,p_0),p_0).
\end{aligned}
\end{equation*}
We observe that the following diagrams commute,
\[
\xymatrix@C=0.5em@R=3em{
	(q_0,p_0) \ar[rr]^{\tilde{F}_{H_d^+}} & & (q_1,p_1) & & & & (q_0,p_0) \ar[rr]^{\tilde{F}_{H_d^-}} & & (q_1,p_1)\\
	& (q_0,p_1) \ar[ul]^{\mathbb{F}^- H_d^+} \ar[ur]_{\mathbb{F}^+ H_d^+} & & & & & & (q_1,p_0) \ar[ul]^{\mathbb{F}^- H_d^-} \ar[ur]_{\mathbb{F}^+ H_d^-}
}
\]

As such, the discrete left and right Hamiltonian maps can be expressed in terms of the discrete fiber derivatives,
\begin{equation} \tilde{F}_{H_d^{\pm}}(q_0,p_0) = \mathbb{F}^+H_d^{\pm} \circ (\mathbb{F}^-H_d^{\pm})^{-1} (q_0,p_0) = (q_1,p_1) ,\end{equation}
and this observation together with Lemma \ref{InvertibilityLemma} ensures that the order of accuracy of the integrator is at least of the order to which the discrete Hamiltonian $H_d^{\pm}$ approximates the exact discrete Hamiltonian $H_d^{\pm,E}$.

However, the Poincar\'e transformed Hamiltonian might be degenerate so we cannot apply the usual implicit function theorem, and need to establish the invertibility of the discrete Legendre transform $\mathbb{F}^-H_d^{\pm}$ in a different way. 

The strongest general result we have been able to establish involves the case where the original Hamiltonian is autonomous, i.e., $H=H(q,p)$, and nondegenerate, and the monitor function is autonomous as well. These assumptions hold for an interesting and useful class of problems, and we will show that the exact discrete left and right Hamiltonians can be reduced to a particular form and that the extended variables $p^t_1$ and $q_1^t$ can be solved for explicitly. As a result, the implicit function theorem is not needed with respect to these variables. 

Hamilton's equations of the Poincar\'e transformed Hamiltonian $\bar{H}(\bar{q},\bar{p})=g(q,p)\left( H(q,p) +p^t \right)$ are given by
\begin{equation*}
	\dot{\bar{q}} = \begin{bmatrix} \nabla_pg(q,p)(H(q,p)+p^t)+ \frac{\partial H}{\partial p}  g(q,p) \\ g(q,p) \end{bmatrix} , \quad \dot{\bar{p}} = -\begin{bmatrix} \nabla_qg(q,p)(H(q,p)+p^t)+\frac{\partial H}{\partial q} g(q,p) \\ 0 \end{bmatrix} .
\end{equation*}

Using these equations, the corresponding exact discrete Hamiltonians are of the form
\begin{equation}
	\bar{H}^{+,E}_d(\bar{q}_0,\bar{p}_1;h) = p_1^\top q_1 + p_1^tq_1^t - \int_0^h \left( p(\tau)^\top \dot{q}(\tau) - g(q(\tau),p(\tau))H(q(\tau),p(\tau)) \right) d\tau ,
\end{equation}
\begin{equation}
	\bar{H}^{-,E}_d(\bar{q}_1,\bar{p}_0;h) = -p_0^\top q_0 - p_0^tq_0^t -  \int_0^h \left( p(\tau)^\top \dot{q}(\tau) - g(q(\tau),p(\tau))H(q(\tau),p(\tau)) \right) d\tau .
\end{equation}

As a result, only one part of these exact discrete left and right Hamiltonians requires approximations of the extended variable $q^t$ and $p^t$. Furthermore, $\dot{p}^t = 0$ implies that $p^t_1=p^t_0$.

Now, let $\bar{H}^{\pm}_d$ be approximations to the exact discrete left and right Hamiltonians of the form
\begin{equation} 
	\begin{aligned} \bar{H}^+_d(\bar{q}_0,\bar{p}_1;h) & = p_1^\top \hat{q}_1(q_0,p_1;h) + p_1^t\hat{q}_1^t(q_0^t,q_0,p_1;h) - I_1(q_0,p_1;h), \\ \bar{H}^-_d(\bar{q}_1,\bar{p}_0;h)  &= - p_0^\top \hat{q}_0(q_1,p_0;h)  - p_0^t\hat{q}_0^t(q_1^t,q_1,p_0;h) - I_2(q_1,p_0;h), \end{aligned} \end{equation}
where $\hat{\cdot}$ denotes an approximation and where $I_1(q_0,p_1;h)$  and $I_2(q_1,p_0;h)$ both approximate
\begin{equation}\int_0^h \left( p(\tau)^\top \dot{q}(\tau) - g(q(\tau),p(\tau))H(q(\tau),p(\tau)) \right) d\tau .\end{equation}

Then, the discrete right Legendre transforms give the following relations for $p_1^t$ and $q_1^t$,
\begin{equation}
\begin{aligned}
	\begin{bmatrix} p_0 \\ p_0^t \end{bmatrix} = \begin{bmatrix} \frac{\partial \hat{q}_1}{\partial q_0}^\top p_1 + p_1^t\frac{\partial \hat{q}_1^t}{\partial q_0} - \frac{\partial I_1}{\partial q_0} \\ \frac{\partial \hat{q}_1^t}{\partial q_0^t}p_1^t \end{bmatrix}, \qquad 
	\begin{bmatrix} q_1 \\ q_1^t \end{bmatrix} = \begin{bmatrix} \hat{q}_1 + \frac{\partial \hat{q}_1}{\partial p_1}^\top p_1 + \frac{\partial \hat{q}_1}{\partial p_1}^\top p_1^t - \frac{\partial I_1}{\partial p_1} \\ \hat{q}_1^t \end{bmatrix}.
\end{aligned}
\end{equation}
Now, since the analytic solution satisfies $p^t_1=p^t_0$, there is no need to approximate $p_1^t$. Therefore, $\frac{\partial \hat{q}_1^t}{\partial q_0^t}=1$. The resulting two systems can both be solved by first setting $p^t_1=p^t_0$, then implicitly solving for $p_1$ in terms of $(q_0^t,q_0,p_1^t,p_1)$, explicitly solving for $q_1$ and finally explicitly solving for $q_1^t$. Since $p_1$ is not determined by $q_1^t$, the implicit function theorem is simply needed for finding $p_1$. Therefore, we need $\operatorname{det}\left(\frac{\partial^2 \bar{H}}{\partial p^2}\right) \neq 0$, and from equation (\ref{d2Hdp2}), this is the same as $\operatorname{det}\left(\frac{\partial H}{\partial p}\nabla_pg(q,p)^\top +g(q,p)\frac{\partial^2 H}{\partial p^2}+\nabla_pg(q,p)\frac{\partial H}{\partial p}^\top \right) \neq 0$. Note this holds for nondegenerate Hamiltonians $H$ and $p$-independent monitor functions. 

Similarly, the discrete left Legendre transforms give the following relations for $p_1^t$ and $q_1^t$,
\begin{equation}
\begin{aligned}
	\begin{bmatrix} p_1 \\ p_1^t \end{bmatrix} = \begin{bmatrix}   \frac{\partial \hat{q}_0}{\partial q_1}^\top p_0 + p_0^t\frac{\partial \hat{q}_0^t}{\partial q_1} + \frac{\partial I_2}{\partial q_1} \\  \frac{\partial \hat{q}_0^t}{\partial q_1^t}p_0^t \end{bmatrix}, \qquad 
	\begin{bmatrix} q_0 \\ q_0^t \end{bmatrix} = \begin{bmatrix} \hat{q}_0 + \frac{\partial \hat{q}_0}{\partial p_0}^\top p_0 + \frac{\partial \hat{q}_0}{\partial p_0}^\top p_0^t + \frac{\partial I_2}{\partial p_0} \\ \hat{q}_0^t \end{bmatrix},
\end{aligned}
\end{equation}
which can be solved provided $\operatorname{det}\left(\frac{\partial H}{\partial p}\nabla_pg(q,p)^\top +g(q,p)\frac{\partial^2 H}{\partial p^2}+\nabla_pg(q,p)\frac{\partial H}{\partial p}^\top \right) \neq 0$. \\

The results that we have established are summarized in the following theorem:

\begin{theorem}
	Given a nondegenerate Hamiltonian $H$, and a monitor function $g \in C^1([0,h])$, such that $\operatorname{det}\left(\frac{\partial H}{\partial p}\nabla_pg(q,p)^\top +g(q,p)\frac{\partial^2 H}{\partial p^2}+\nabla_pg(q,p)\frac{\partial H}{\partial p}^\top \right)\neq 0$. Then, if the discrete Hamiltonian $\bar{H}^{\pm}_d$, approximates the exact discrete Hamiltonian $\bar{H}_d^{\pm,E}$ to order $r$, i.e.,
	\begin{equation} 
		\bar{H}^{\pm}_d(\bar{q}_0, \bar{p}_1;h) =\bar{H}_d^{\pm,E}(\bar{q}_0,\bar{p}_1;h)+\mathcal{O}(h^{r+1}) ,
	\end{equation}
	then the discrete Hamiltonian map $\tilde{F}_{\bar{H}^{\pm}_d}:(\bar{q}_k,\bar{p}_k)\mapsto(\bar{q}_{k+1},\bar{p}_{k+1})$, viewed as one-step method, is order $r$ accurate. 
\end{theorem}

\begin{remark}
	It should be noted that the assumptions that the original Hamiltonian is nondegenerate and autonomous fail to hold in the application we consider of time-adaptive variational integrators to the discretization of the Bregman Hamiltonian associated with accelerated optimization, as it is time-dependent. This is unavoidable, as it models a system with dissipation, which cannot be described with an autonomous Hamiltonian as the Hamiltonian would otherwise be an integral of motion as it is the Noether quantity associated with time translational symmetry. In the cases when the original Hamiltonian is degenerate or nonautonomous, we need to analyze the solvability of the discrete Hamiltonian equations on a case-by-case basis, but as we demonstrate, this can be done in the case of the Bregman Hamiltonian with the given choices of monitor function $g(t)$ and discrete Hamiltonians that we consider.
\end{remark}

\subsection{Numerical Tests on Kepler's Planar 2-Body Problem}\label{KeplerSection}
We will now demonstrate the approach using Hamiltonian Taylor variational integrators, presented in Section \ref{HTVISection}, on Kepler's planar 2-Body Problem. For a lucid exposition, we will at first assume that $g(q,p) = g(q)$ and $H(q,p)=\frac{1}{2}p^\top M^{-1}p+V(q)$. Consider the discrete right Hamiltonian given by approximating $\bar{q}_1$ with a first-order Taylor method about $\bar{q}_0$, approximating $\bar{p}_0$ with a zeroth-order Taylor expansion about $\bar{p}_0$, and using the rectangular quadrature rule about the initial point:
\begin{equation}
	\bar{H}_d^+ = p_1^\top \left(q_0+\frac{1}{2} hg(q_0)M^{-1}p_1\right) + p_1^t(q_0^t+hg(q_0)) + hg(q_0) V(q_0)  .
\end{equation} 

The corresponding variational integrator is given by
\begin{equation} \label{adeulb}
	\begin{aligned}
		\bar{p}_1 &= \begin{bmatrix} p_0-hg(q_0)\nabla V(q_0)-h\nabla g(q_0)\left( \frac{1}{2} p_1^\top M^{-1}p_1+ V(q_0)+p_0^t \right) \\ p_0^t \end{bmatrix}, \\
		\bar{q}_1 &= \begin{bmatrix} q_0+hg(q_0)M^{-1}p_1 \\ q_0^t + hg(q_0) \end{bmatrix} .
	\end{aligned} 
\end{equation} 

This integrator is merely symplectic Euler-B applied to the transformed Hamiltonian system
\begin{equation}
	\bar{q}_1  = \bar{q}_0 + h\frac{\partial \bar{H}(\bar{q}_0,\bar{p}_1)}{\partial \bar{p}} , \qquad 
	\bar{p}_1  = \bar{p}_0 - h\frac{\partial \bar{H}(\bar{q}_0,\bar{p}_1)}{\partial \bar{q}} .
\end{equation}

In fact, this is precisely the adaptive symplectic integrator first proposed in \cite{Ha1997} and presented in \cite[page 254]{LeRe2005}. Most existing symplectic integrators can be interpreted as variational integrators, but there are also new methods that are most naturally derived as variational integrators. We will also consider a fourth-order Hamiltonian Taylor variational integrator (HTVI4), which is distinct from any existing symplectic method. 

One of the most important aspects of implementing a variable time-step symplectic integrator of this form is a well chosen monitor function, $g(q)$. We need $g$ to be positive-definite, so that we never stall or march backward in time. Noting that the above integrator is first-order, a natural choice is to use the second-order truncation error given by $-\frac{(q_1^t-q_0^t)^2}{2}M^{-1}\nabla V(q_0)$. Let $tol$ be some desired level of accuracy. Then, using $q_1^t-q_0^t = hg(q_0)$, one choice for $g$ would be,
\begin{equation}
	g(q_0) = \frac{tol}{\| \frac{(q_1^t-q_0^t)^2}{2}g(q_0)M^{-1}\nabla V(q_0) \|} = \frac{tol}{\| \frac{h^2g(q_0)^3}{2}M^{-1}\nabla V(q_0) \|} ,
\end{equation}
Thus 
\begin{equation}
	g(q_0) = \left(\frac{tol}{\| \frac{h^2}{2}M^{-1}\nabla V(q_0) \|} \right)^{\frac{1}{4}}.
\end{equation}

Experimentally, the $4$th-root did not affect results very much, but required messier computations, which is the reason why we have chosen the simpler yet very similar monitor function 
\begin{equation}
	g(q_0) = \frac{tol}{\| \frac{h^2}{2}M^{-1}\nabla V(q_0) \|}, \label{trunerror}
\end{equation}
which achieves an error which is comparable to the chosen value of $tol$. 

Alternative choices for the monitor function $g(q)$, proposed in \cite{Ha1997}, include the $p$-independent arclength parameterization
\begin{equation}
	g(q) = (2(H_0 - V(q)) + \nabla V(q)^\top M^{-1}\nabla V(q))^{-\frac{1}{2}}, \label{arclength}
\end{equation}
and a choice particular to Kepler's two-body problem,
\begin{equation}
	g(q) = q^\top q, \label{qnorm}
\end{equation}
which is motivated by Kepler's second law, which states that a line segment joining the two bodies sweeps out equal areas during equal intervals of time.

\begin{figure}[!h]
	\hspace*{-25mm}
	\centering
	\begin{minipage}[b]{0.55\textwidth}
		\includegraphics[width=\textwidth]{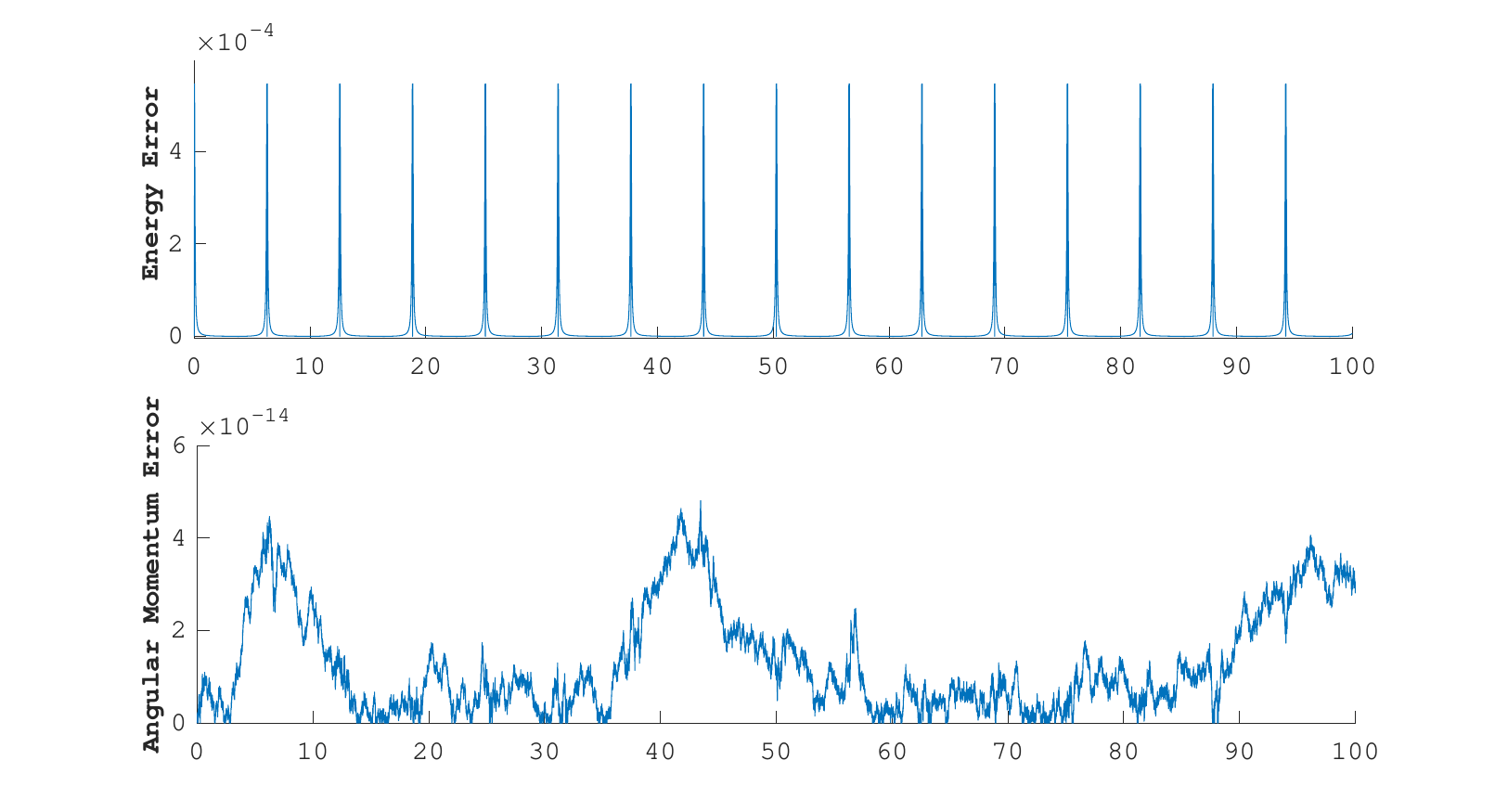}
		\subcaption{	Symplectic Euler-B}
	\end{minipage}
	\hspace{-10mm}
	\begin{minipage}[b]{0.55\textwidth}
		\includegraphics[width=\textwidth]{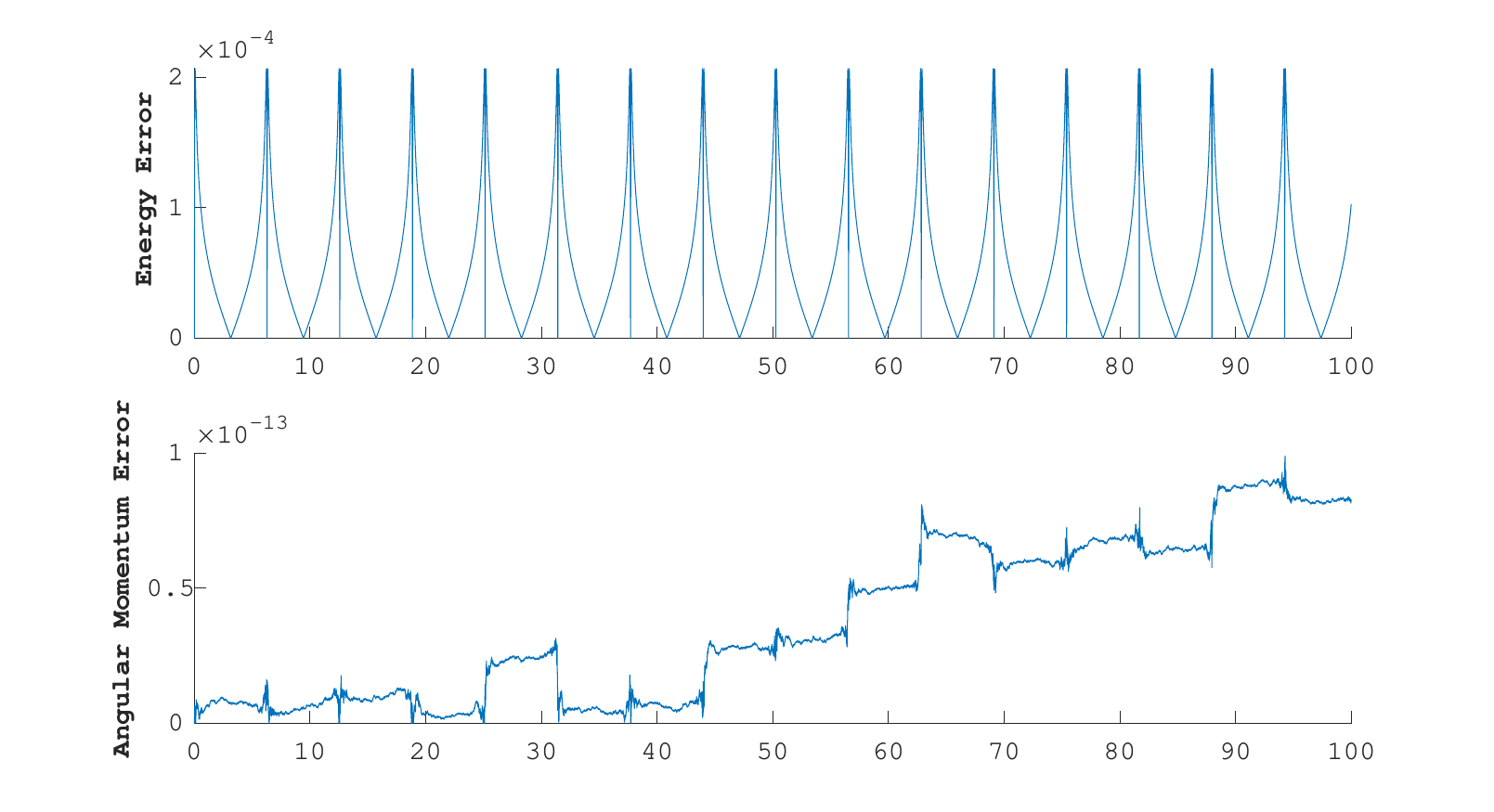}
		\subcaption{Adaptive method with monitor function \eqref{trunerror} }
	\end{minipage}
	\hspace*{-25mm}
	\caption{Symplectic Euler-B and the adaptive algorithm with monitor function given by equation (\ref{trunerror}) were applied to Kepler's planar two-body problem over a time interval of $[0,100]$ with an eccentricity of 0.9.\label{Kepler1}}
\end{figure}

\begin{figure}[!h]
	\includegraphics[width=1\textwidth]{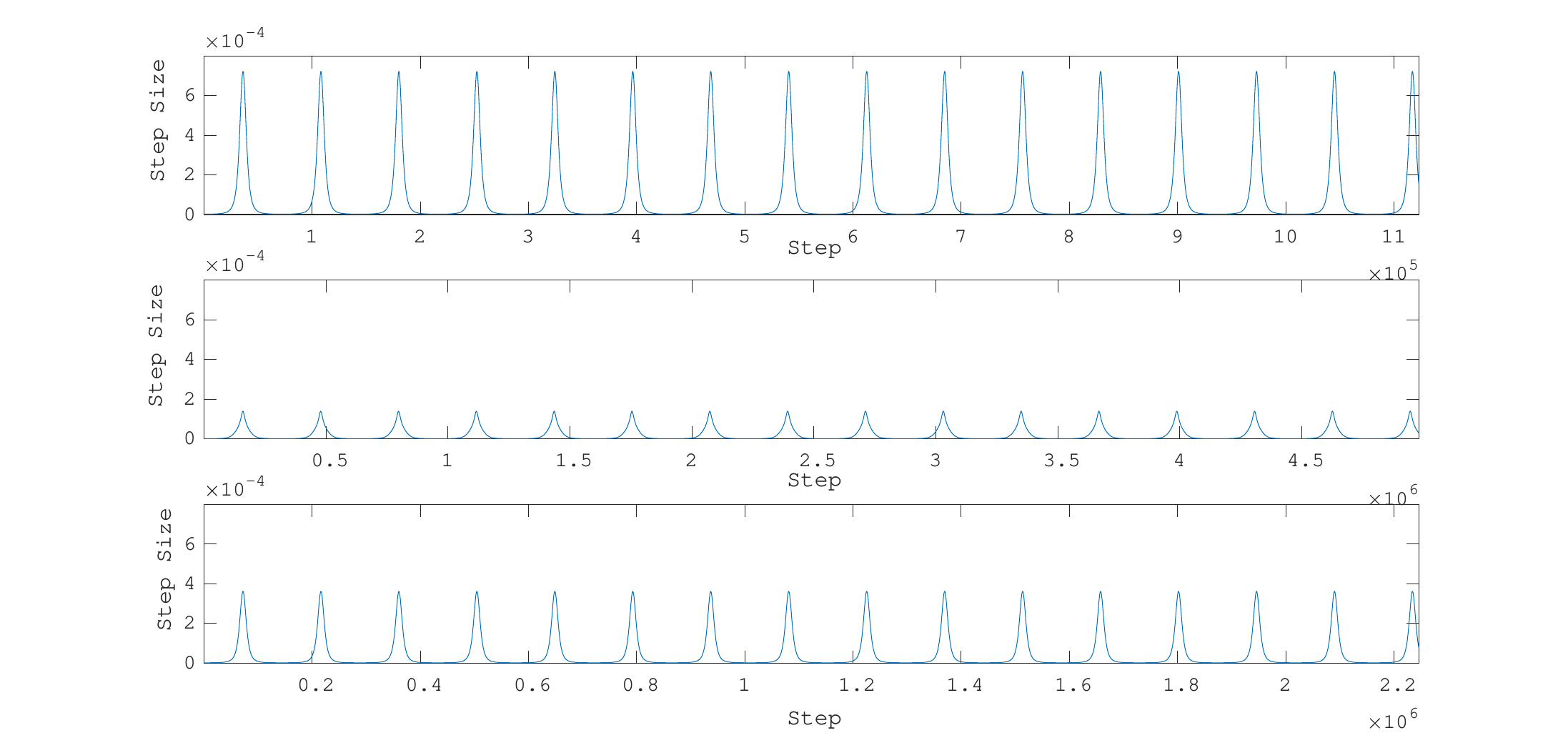}
	\caption
	{Time-steps taken for the various choices of monitor functions. The top, middle and bottom plots correspond to the monitor functions \eqref{trunerror}, \eqref{arclength}, and \eqref{qnorm}, respectively. All of the monitor functions appear to increase and decrease the time-step at the same points along the trajectory, but clearly \eqref{trunerror} allowed for the larger steps to be taken.\label{EBstep}}
\end{figure}

We have tested the algorithm given by equation \eqref{adeulb} on Kepler's planar two-body problem, with an eccentricity of 0.9, using the three choices of monitor function $g$ given by \eqref{trunerror}, \eqref{arclength}, and \eqref{qnorm}. Of these three choices, \eqref{qnorm} is specific to Kepler's two-body problem, while \eqref{trunerror} and \eqref{arclength} are more general choices. Unlike \eqref{arclength} which is independent of the order of the method, \eqref{trunerror} is based on the truncation error and thus the corresponding cost of computing this function will increase as the order of the method increases. Simulations using Kepler's two-body problem with an eccentricity of 0.9 over a time interval of $[0,1000]$ were run using the three different choices of $g$ and the usual symplectic Euler-B. Results indicate that symplectic Euler-B takes the most steps and computational time to achieve a level of accuracy around $10^{-5}$. To achieve this level of accuracy, the choice of the truncation error monitor function, \eqref{trunerror}, resulted in the least number of steps, and the second lowest computational time. The lowest computational time belonged to \eqref{qnorm}, but it used significantly more steps than \eqref{trunerror}. The lower computational cost can be attributed to the cheaper evaluation cost of the monitor function and its derivative. Finally, the monitor function \eqref{arclength} required the most steps and computational time of the adaptive algorithms, but it is still a good choice in general given its broad applicability. Figures \ref{Kepler1} and \ref{EBstep} present the energy and angular momentum errors for the fixed timestep method vs. adaptive timestep method, and the time-steps for the different monitor functions, respectively.


Next, we consider a Type II fourth-order Hamiltonian Taylor variational integrator (HTVI4) constructed using the strategy from Section \ref{HTVISection} and the automatic differentiation package from \cite{Neid2010,10.1145/2450153.2450155}. We will drop the assumption of $p$-independent monitor functions and consider $g(q,p)$. The following monitor functions were considered,
\begin{align}
	g(q) &= \left( q^\top q \right)^{\gamma} \ \text{for} \ \gamma = 0.5 \text{ and }1 , & (Gamma) \label{gamma} \\
	g(q) &= \left( 2(H_0-V(q))+\nabla V(q)^\top M^{-1}\nabla V(q)\right)^{-\frac{1}{2}} ,  & (Arclength) \label{arc_length} \\ 
	g(q,p) &= \| p^t-L(q,M^{-1}p) \|_2^{-1} . & (Energy) \label{energy}
\end{align} 
The monitor function \eqref{energy} was originally intended to be $\|p^t+H(q,p)\|_2^{-1}$, but experimental results suggested that \eqref{energy} is the better choice. We will discuss the shortcomings of using the inverse energy error in the next paragraph. Note that $\| L(q,M^{-1}p) \|_2^{-1}$ also performs decently, but the addition of $p^t=-H(q_0,p_0)$ showed noticeable improvement. It was noted in \cite{Ha1997} that the inverse Lagrangian has been considered as a possible choice for $g$ in the Poincar\'e transformation, but not in the framework of symplectic integration. While the choice of \eqref{gamma} was generally the most efficient, \eqref{energy} was very close in terms of efficiency and offers a more general monitor function. This also implies that efficiency is not limited to only $q$ or $p$-independent monitor functions. However, various attempts to construct separable transformed Hamiltonians (see \cite{BlBu2004, BlIs2012}) required the use of $q$ or $p$-independent monitor functions, so this is where such monitor functions are most useful. 

In the case of monitor functions involving the gradient, higher-order derivatives will be required for higher-order Taylor variational integrators, but there are efficiencies to be had when leveraging the higher-order derivatives already being calculated for the underlying Taylor method and Hessian-vector multiplication that can be done efficiently without needing to explicitly construct the full Hessian \cite{10.1093/imanum/12.2.135}. The calculation of higher-order derivatives do come with a higher cost, and in the case of Kepler's 2-body problem there is a clear computational advantage in using the gradient-free gamma monitor function \eqref{gamma} as shown in Tables \ref{table1} and \ref{table2}. However, the gamma monitor function \eqref{gamma} is more specific to Kepler's 2-body, while the energy and arclength monitor functions are applicable to a wider range of problems. Monitor functions that are both general and efficient would be highly desirable. 

\begin{figure}[!h]
	\includegraphics[width=0.92\textwidth]{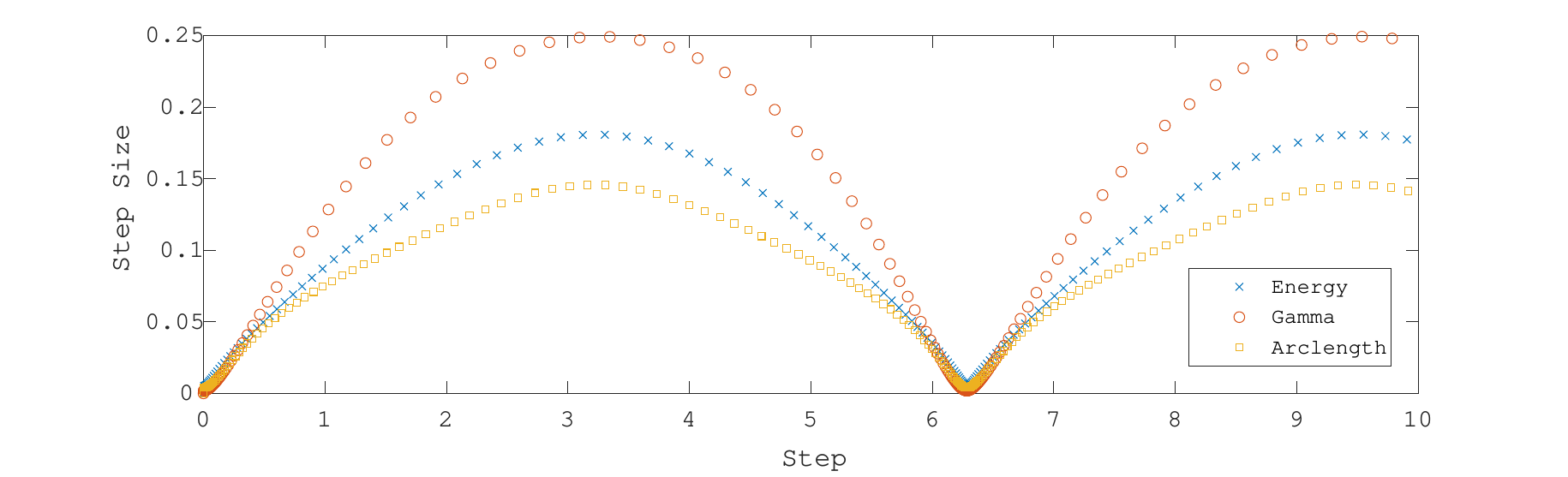}
	\caption
	{Time-steps taken for the various choices of monitor functions. The energy \eqref{energy} and gamma \eqref{gamma} monitor functions  performed better, in terms of fewest steps, lowest computational cost and lowest global error, than the arclength monitor function \eqref{arc_length}. Note \eqref{energy} did not take the largest nor the smallest steps.}
	\label{Steps9}
\end{figure}

Figure \ref{Steps9} displays the time-steps taken for the different choices of monitor functions for this fourth-order Hamiltonian Taylor variational integrator.

The truncation error monitor function, \eqref{trunerror}, performed quite well for first-order methods, and this motivated the choice of using Taylor variational integrators, since derivatives would be readily available. However, its success cannot as easily be applied to higher-order methods. This is due to the fact that for higher-order truncation errors, one obtains an implicit differential-algebraic definition of the monitor function. This deviates from the first-order case, where the monitor function can be solved for explicitly. Another seemingly natural choice for the monitor function is the inverse of the energy error. However, Taylor variational integrators are constructed using Taylor expansions about the initial point, and consequently the monitor function is mostly evaluated at or near the initial point. If the initial point is at a particularly tricky part of the dynamics and requires a small first step, then the energy error at the first step will not reflect this, since initially the energy error is zero. In contrast, the inverse Lagrangian will be small at an initial point that requires a small first step. The inverse energy error may work well for methods that primarily evaluate the energy error at the end point rather than the initial point. It is also often advantageous to bound the time-step below or above. As noted in \cite[page 248]{LeRe2005}, this can be done by setting $a=\frac{\Delta t_{\text{min}}}{\Delta \tau}$ and $b=\frac{\Delta t_{\text{max}}}{\Delta \tau}$, and defining the new monitor function as $\hat{g} = b\frac{g + a}{g+b}$. 

Tables \ref{table1} and \ref{table2} display a comparison of bounds, computational time, steps, and error. 
We can note that for methods such as the Taylor variational integrator, bounding $g(q,p)$ bounds the time-step, but not directly. Also, compared to non-adaptive variational integrators, such as the non-adaptive Taylor variational integrator and the St\"ormer--Verlet method (SV), the adaptive methods showed a significant gain in efficiency for Kepler's 2-body planar problem with high eccentricity, while low eccentricity models do not need nor do they benefit from adaptivity. A Hamiltonian dynamical system with regions of high curvature in the vector field and its norm will in general benefit from an adaptive scheme such as the one outlined here. 

\begin{table}[ht]
	\centering
	\resizebox{\textwidth}{!}{\begin{tabular}{ccccccccccc}
			\hline
			Method & Monitor $g(q,p)$ & $h$ & min Step& max Step& min $g$ & max $g$ & Energy Error & Global Error & Steps & Time\\
			\hline
			HTVI4 & Gamma & 0.1 & 0.0020 &	0.2493 & 0.01	& 8	& 1.43E-05 & 7.09E-06 & 181	& 26.9\\
			HTVI4 & Energy & 0.1 & 0.0051 &	0.1809 &	1E-04 &	2	& 1.93E-06	& 4.76E-06	& 146 &	28.3\\
			HTVI4 & Arclength & 0.1 & 0.0040 &	0.1458 &	3E-03	& 0.3	& 1.10E-04 & 3.69E-05	& 185	& 70.2\\
			HTVI4 & - & 0.0025 &  0.0025 & 0.0025 & - & - & 2.50E-06	& 2.89E-05 & 4000 & 120\\
			SV & - & 5E-05 &  5E-05 & 5E-05 & - & - & 3.12E-06 & 4.68E-05	& 2E05 & 1.9\\
			\hline
	\end{tabular}}
	\caption{Kepler Planar 2-Body Problem, Eccentricity = 0.9} 
	\label{table1}
\end{table}

\begin{table}[ht]
	\centering
	\resizebox{\textwidth}{!}{\begin{tabular}{ccccccccccc}
			\hline
			Method & Monitor $g(q,p)$ & $h$ & min Step& max Step& min $g$ & max $g$ & Energy Error & Global Error & Steps & Time\\
			\hline
			HTVI4 & Gamma & 0.1 & 6E-05	& 0.2648 & 5E-04	& 8	& 4.88E-05 & 5.60E-06	& 372	& 49.3\\
			HTVI4 & Energy & 0.03 & 1.5E-04 & 0.1462 &	1E-06 & 5 & 9.13E-06 & 4.63E-06 & 383 & 58.4\\
			HTVI4 & Arclength & 0.1 & 5E-05 & 0.1379 &	8E-04 & 10	& 1.31E-05 & 1.49E-05	& 691 & 146.0\\
			HTVI4 & - &  5E-04 &  5E-04 & 5E-04 & - & - & 1.38E-01	& 7.83E-01 & 2E04 & 525.2\\
			SV & - & 5E-07 &  5E-07 & 5E-07 & - & - & 3.34E-06 & 2.68E-05	& 2E07 & 189.2\\
			\hline
	\end{tabular}}
	\caption{Kepler Planar 2-Body Problem, Eccentricity = 0.99} 
	\label{table2}
\end{table}

\section{Application to Symplectic Accelerated Optimization} \label{OptimizationSection}

\subsection{Accelerated Optimization}

Efficient optimization has become one of the major concerns in data analysis. Many machine learning algorithms are designed around the minimization of a loss function or the maximization of a likelihood function. Due to the ever-growing scale of the data sets and size of the problems, there has been a lot of focus on first-order optimization algorithms because of their low cost per iteration. The first gradient descent algorithm was proposed in \cite{Cauchy1847}  by Cauchy to deal with the very large systems of equations he was facing when trying to simulate orbits of celestial bodies, and many gradient-based optimization methods have been proposed since Cauchy's work in 1847. In 1983, Nesterov's Accelerated Gradient method was introduced in \cite{Nes83},
\begin{equation}\label{NesterovUpdate}
	x_k = y_{k-1} -h\nabla f(y_{k-1}) , \qquad y_k = x_k + \frac{k-1}{k+2} (x_k - x_{k-1}), 
\end{equation}
which converges in $\mathcal{O}(1/k^2)$ to the minimum of the convex objective function $f$, improving on the $\mathcal{O}(1/k)$ convergence rate exhibited by the standard gradient descent methods.
This $\mathcal{O}(1/k^2)$ convergence rate was shown in \cite{Nes04} to be optimal among first-order methods using only information about $\nabla f$ at consecutive iterates. This phenomenon in which an algorithm displays this improved rate of convergence is referred to as acceleration, and other accelerated algorithms have been derived since Nesterov's algorithm, such as accelerated mirror descent \cite{Nem1983}, and accelerated cubic-regularized Newton's method \cite{Nes08}. More recently, it was shown in \cite{SuBoCa16} that Nesterov's Accelerated Gradient method limits to the second order ODE, 
\begin{equation}   \ddot{x}(t) + \frac{3}{t} \dot{x}(t) +\nabla f(x(t)) = 0 ,\end{equation} 
as $h\rightarrow 0$. The authors also proved that the objective function $f(x(t))$ converges to its optimal value at a rate of $\mathcal{O}(1/t^2)$ along the trajectories of this ODE. It was then shown in \cite{WiWiJo16} that in continuous time, the convergence rate of $f(x(t))$ can be accelerated to an arbitrary convergence rate $\mathcal{O}(1/t^p)$, by considering flow maps generated by time-dependent Lagrangian and Hamiltonian systems. We will present this result in more detail in the next section together with the variational framework introduced in \cite{WiWiJo16} for accelerated optimization, which will be at the heart of our approach. \\

\subsection{Variational Framework for Accelerated Optimization}

In this section, we will review the variational framework introduced in \cite{WiWiJo16} for accelerated optimization which will be the basis for the methods we will design. In a general space $\mathcal{X}$, given a convex, continuously differentiable function $h:\mathcal{X} \rightarrow \mathbb{R}$ such that $\Vert \nabla h(x) \Vert \rightarrow \infty $ as $\Vert x \Vert \rightarrow \infty $, its corresponding Bregman divergence is given by
\begin{equation}  D_h(x,y) = h(y)-h(x) - \langle \nabla h(x), y-x \rangle .  \end{equation} 
We then define the Bregman Lagrangian and Hamiltonian
\begin{equation} \label{BregmanLGeneral}
	\mathcal{L}_{\alpha,\beta,\gamma}(x,v,t) = e^{\alpha(t) + \gamma(t)} \left[   D_h(x+e^{-\alpha(t)}v , x) - e^{\beta(t)} f(x)  \right]  ,
\end{equation}
\begin{equation} \label{BregmanHGeneral}
	\mathcal{H}_{\alpha,\beta,\gamma}(x,r,t) = e^{\alpha(t) + \gamma(t)} \left[   D_{h^*}(\nabla h(x)+e^{-\gamma(t)}r , \nabla h(x) ) + e^{\beta(t)} f(x)  \right] ,
\end{equation}
which are scalar valued functions of position $x\in \mathcal{X}$, velocity $v\in \mathbb{R}^d$,  momentum $r\in \mathbb{R}^d$, and time $t$, which are parametrized by smooth functions of time, $\alpha,\beta,\gamma$, and where $h^* = \sup_{v\in T\mathcal{X}}{\left[ \langle r, v \rangle - h(v) \right]}$ is the Legendre transform (or convex dual function) of $h$. These parameters $\alpha,\beta,\gamma$ are said to satisfy the ideal scaling conditions if
\begin{equation}\label{IdealScaling}
	\dot{\beta}(t) \leq e^{\alpha (t)} \qquad \text{and} \qquad \dot{\gamma}(t) = e^{\alpha (t)}.
\end{equation}
If the ideal scaling conditions are satisfied, then by Theorem 1.1 in \cite{WiWiJo16},
\begin{equation}
	f(x(t))-f(x^*) \leq \mathcal{O}(e^{-\beta(t)}).
\end{equation}
Another very important property of this family of Bregman Lagrangians is its closure under time dilation, proven in Theorem 1.2 of \cite{WiWiJo16}:
\begin{theorem}\label{ThmTimeDilation}
	If $x(t)$ satisfies the Euler-Lagrange equations corresponding to the Bregman Lagrangian $\mathcal{L}_{\alpha,\beta,\gamma}$, then the reparametrized curve $y(t) = x(\tau(t))$  satisfies the Euler-Lagrange equations corresponding to the Bregman Lagrangian $\mathcal{L}_{\tilde{\alpha},\tilde{\beta},\tilde{\gamma}}$ where 
	\begin{equation}\label{TimeDilationTransformation}
		\tilde{\alpha}(t) = \alpha(\tau(t)) + \log{\dot{\tau} (t)},  \qquad \tilde{\beta}(t) = \beta(\tau(t)), \qquad \tilde{\gamma}(t) = \gamma(\tau(t)) ,
	\end{equation}
	and
	\begin{equation} \mathcal{L}_{\tilde{\alpha},\tilde{\beta},\tilde{\gamma}} (x,v,t) = \dot{\tau} (t) \mathcal{L}_{\alpha,\beta,\gamma}\left(x,\frac{1}{\dot{\tau} (t)}v,\tau(t)\right) .  \end{equation}	
	Furthermore $\alpha,\beta,\gamma$ satisfy the ideal scaling equation (\ref{IdealScaling}) if and only if $\tilde{\alpha},\tilde{\beta},\tilde{\gamma}$ do.
\end{theorem} 

A subfamily of Bregman Lagrangians of interest, indexed by a parameter $p>0$, is given by the choice of functions
\begin{equation} \label{AlphaBetaGamma}
	\alpha(t) = \log{p} - \log{t} , \qquad  \beta(t) = p \log{t} + \log{C},  \qquad \gamma(t) = p \log{t},
\end{equation}
where $C>0$ is a constant. The Bregman Lagrangian and Hamiltonian become
\begin{equation} 
	\mathcal{L}(x,v,t) = pt^{p-1} \left[   D_{h}\left(x+\frac{t}{p} v,x\right) - Ct^p f(x)  \right] ,
\end{equation}  
\begin{equation} 
	\mathcal{H}(x,r,t)= pt^{p-1} \left[   D_{h^*}(\nabla h(x)+t^p r , \nabla h(x) ) + Ct^pf(x)  \right] .
\end{equation}
These parameter functions are of interest since they satisfy the ideal scaling equation (\ref{IdealScaling}), and the resulting evolution $x(t)$ of the corresponding dynamical system was shown in \cite{WiWiJo16} to satisfy the aforementioned  $\mathcal{O}(1/t^p) $ convergence rate, \begin{equation}  f(x(t))-f(x^*) \leq \mathcal{O}(1/t^p), \end{equation}  where $x^*$ is the desired minimizer of the objective function $f$.  \\

\subsection{Adaptive Variational Integrators for Symplectic Accelerated Optimization}

For simplicity of exposition, we will consider the case where $h(x) = \frac{1}{2} \langle x,x \rangle $. Our new approaches will make use of the adaptive framework developed in Section \ref{AdaptiveSection} via carefully chosen Poincar\'e transformations. Recalling the discussion of Section \ref{PoincareSubsection}, there might not be a Lagrangian formulation for the future Poincar\'e transformed systems, so we will need to work from the Hamiltonian point of view to design variational integrators. When $h(x) = \frac{1}{2}  \langle x,x \rangle $, the Bregman Hamiltonian with parameters $\alpha,\beta,\gamma$ given by equation (\ref{AlphaBetaGamma}) for a specific value of $p>0$ becomes
\begin{equation} \label{BregmanH}
	\mathcal{H}(q,r,t)= \frac{p}{2t^{p+1}} \langle r , r \rangle + Cpt^{2p-1} f(q).
\end{equation}  
As mentioned in the previous section, the solution to the corresponding Hamilton's equations was shown in \cite{WiWiJo16} to satisfy the convergence rate  \begin{equation}  f(q(t))-f(q^*) \leq \mathcal{O}(1/t^p) ,\end{equation}  where $q^*$ is the desired minimizer of the objective function $f$. Together with the time-dilation result from Theorem \ref{ThmTimeDilation}, this implies that this entire subfamily of Bregman trajectories indexed by the parameter $p$ can be obtained by speeding up or slowing down along the Bregman curve in spacetime corresponding to any specific value of $p$. We will present two new approaches based on the adaptive framework developed in Section \ref{AdaptiveSection}  to integrate the Bregman Hamiltonian dynamics, thereby solving the optimization problem, and then compare their performance. \\

\noindent \underline{\textbf{Direct Approach}}: Our first approach will use our adaptive framework with monitor function $g(q,r) = 1$ to design a variational integrator for the Bregman Hamiltonian given in equation (\ref{BregmanH}) for a given value of $p>0$. This choice of monitor function will convert the time-dependent Bregman Hamiltonian into an autonomous Hamiltonian in extended phase space. More precisely, given a value of $p>0$, the time transformation $t \mapsto \tau$ given by $\frac{dt}{d\tau} = g(q,t,r) = 1$ generates the Poincar\'e transformed Hamiltonian
\begin{equation}\label{pHamiltonian}
	\bar{H}(\bar{q},\bar{r}) = \frac{p}{2(q^t)^{p+1}} \langle r , r \rangle + Cp(q^t)^{2p-1} f(q) + r^t ,
\end{equation}
in the phase space with extended coordinates $(\bar{q},\bar{r})$. This strategy is equivalent to the usual trick to remove time-dependency by considering time $t$ as an additional position variable and adding a corresponding conjugate momentum variable, which is the energy (see \cite{JordanSymplecticOptimization} for an example with Hamiltonian given by equation (\ref{BregmanHGeneral})). This shows that our adaptive framework is very general and can also be used for other purposes than solely enforcing a desired variable time stepping.\\

\noindent\underline{\textbf{Adaptive Approach}}: Our second approach will exploit the time-dilation property of the Bregman dynamics together with our adaptive framework with a carefully tuned monitor function. More precisely, we will use adaptivity to transform the Bregman Hamiltonian corresponding to a specific value of $p>0$ into an autonomous version of the Bregman Hamiltonian corresponding to a smaller value $\mathring{p} < p$ in extended phase-space. This will allow us to integrate the higher-order Bregman dynamics corresponding to the value $p$ while benefiting from the computational efficiency of integrating the lower-order Bregman dynamics corresponding to the value $\mathring{p}<p$. Explicitly, solving equation (\ref{TimeDilationTransformation}) for $\tau(t)$ to transform the Bregman dynamics corresponding to the values of $\alpha,\beta,\gamma$ as in equation (\ref{AlphaBetaGamma}) for a given value of $p$ into the Bregman dynamics corresponding to the values of $\tilde{\alpha},\tilde{\beta},\tilde{\gamma}$ as in equation (\ref{AlphaBetaGamma}) for a given value $\mathring{p}<p$ yields $\tau(t) = t^{\mathring{p}/p}$. The corresponding monitor function is given by 
\begin{equation} 	\frac{dt}{d\tau} = g(q,t,r) = \frac{p}{\mathring{p}} t^{1-\frac{\mathring{p}}{p}},\end{equation} 
and generates the Poincar\'e transformed Hamiltonian
\begin{equation} \label{ADHamiltonian}
	\bar{H}(\bar{q},\bar{r}) =  \frac{1}{\mathring{p}} \left[ \frac{p^2}{2(q^t)^{p+\frac{\mathring{p}}{p}} }  \langle r , r \rangle  + Cp^2 (q^t)^{2p-\frac{\mathring{p}}{p}} f(q) + p r^t (q^t)^{1-\frac{\mathring{p}}{p}}    \right] . 
\end{equation}

\hfill

\subsection{Presentation of Numerical Methods} \label{SectionMethods}

In this section, we will test the performance of our new adaptive framework by implementing variational and non-variational integrators in the case where $\mathcal{X} = \mathbb{R}^d$ and $ \langle x,x\rangle = x^\top x$, and discuss the results obtained. Keeping in mind the machine learning application where data sets are very large, we will restrict ourselves to explicit first-order optimization algorithms. 

Looking at the forms of Hamilton's equations in both the Direct and Adaptive approach, we can note that the objective function $f$ and its gradient $\nabla f$ only appear in the expression for $\dot{\bar{r}}$, and are functions of $q$ only. Looking back at the construction of HTVIs from Section \ref{HTVISection}, denoting the order of the Taylor methods by $\rho$ instead of $r$ to avoid confusion with the extended momentum variable $\bar{r}$, we can note that both the Type II and Type III approaches require $\rho$-order Taylor approximations of $r$ and $(\rho +1)$-order Taylor approximations of $q$. This means that the highest value of $\rho$ that we can choose to obtain a gradient-based algorithm is $\rho = 1$.  Now, the starting point of the Type III approach is a $(\rho+1)$-order Taylor approximation $\tilde{q}_0$ of $q_0$ around $q_0$. As a consequence, the subsequent steps in the Type III method with $\rho =0$ and $\rho = 1$ will contain evaluations of the objective function $f$ and its gradient $\nabla f$ at this approximation $\tilde{q}_0$. Aside from the inconvenience of the function evaluations not being at the iterates $q_0$ and $q_1$ themselves, if $f$ is a nonlinear function, this will also generate nonlinearity in the equations for the updates. As a result, we will not be able to design an explicit algorithm, or at least not a general explicit algorithm that would work for all the functions $f$ considered. On the other hand, the starting point of the Type II approach is a $\rho$-order Taylor approximation $\tilde{p}_0$ of $p_0$ around $p_0$. A similar issue as for the Type III case arises when $\rho = 1$ due to the approximations $(q_{c_i},p_{c_i})  = \Psi_{c_i h}^{(\rho)}(\tilde{q}_0,p_0) $. Therefore, we cannot design a general explicit algorithm for the Type II case with $\rho = 1$. The remaining possibility is to construct a Type II HTVI using $\rho = 0$. This will produce explicit gradient-based algorithms, where all the evaluations of the objective function $f$ and its gradient $\nabla f$ are performed at the iterates $q_0$ and $q_1$.  Note that when $\rho = 0 $, we have $(q_{c_i},p_{c_i})  = \Psi_{c_i h}^{(0)}(\tilde{q}_0,p_0) = (\tilde{q}_0,p_0) $ for all $i$, so for given values of $p$ and $\mathring{p}$, every quadrature rule generates the same integrator. Following the method of Section \ref{HTVISection}, with time-step $h$, we will derive explicit gradient-based HTVIs: \\

\noindent \underline{\textbf{Type II Hamiltonian Taylor Variational Integrators (HTVIs) with $\rho = 0$}}: As mentioned earlier, since $\rho = 0$, the choice of quadrature rule does not matter, so we can take the rectangular quadrature rule about the initial point  ($c_1 = 0$ and $b_1 =1$). We approximate $\bar{r}(0)=\tilde{\bar{r}}_0$ via 
$ \bar{r}_1 = \pi_{T^*\bar{Q}} \circ \Psi_h^{(0)}(q_0,\tilde{\bar{r}}_0) = \tilde{\bar{r}}_0$ and generate approximations $	(\bar{q}_{c_1},\bar{r}_{c_1})   = \Psi_{c_1 h}^{(0)}(\bar{q}_0,\tilde{\bar{r}}_0) = (\bar{q}_0,\tilde{\bar{r}}_0) $. 

%\nolinenumbers  

\begin{table}[!h]
	\resizebox{15cm}{!} {
	\begin{tabular}{l|l}
		\textbf{Direct Approach} & \textbf{Adaptive Approach}  \\ \hline
		$	\tilde{\bar{q}}_1  = \pi_ {\bar{Q}}  \circ \Psi_h^{(1)}(\bar{q}_0,\tilde{\bar{r}}_0) =  \begin{bmatrix}   q_0 + h \frac{p r_0}{(q_0^t)^{p+1}}  \\ q_0^t + h \end{bmatrix} $,  & $  \tilde{\bar{q}}_1  = \pi_ {\bar{Q}}  \circ \Psi_h^{(1)}(\bar{q}_0,\tilde{\bar{r}}_0) = \begin{bmatrix}   q_0 + h \frac{p^2 r_0}{\mathring{p}} (q_0^t)^{-p-\frac{\mathring{p}}{p}} \\ q_0^t + h \frac{p}{\mathring{p}} (q_0^t)^{1-\frac{\mathring{p}}{p}} \end{bmatrix} $	,
		\\
		$ \begin{aligned} H_d^+(\bar{q}_0,\bar{r}_1;h) &= r_1^\top  q_0 + r_1^t q_0^t + h  \frac{p}{2(q_0^t)^{p+1}} r_1^\top r_1 \\ &  \qquad + hCp(q_0^t)^{2p-1} f(q_0)  + hr_1^t . \end{aligned} $
		
		&       \small $ \begin{aligned}  & H_d^+(\bar{q}_0,\bar{r}_1;h) = r_1^\top  q_0 + r_1^t q_0^t + h \frac{p^2}{2\mathring{p} (q_0^t)^{p+\frac{\mathring{p}}{p}}  }   r_1^\top r_1 \\ & \qquad   \qquad \qquad + hC \frac{p^2}{\mathring{p}} (q_0^t)^{2p-\frac{\mathring{p}}{p}} f(q_0) +   h \frac{p}{\mathring{p}}  (q_0^t)^{1-\frac{\mathring{p}}{p}}  r_1^t  . \end{aligned}$ \\
		
		The discrete right Hamilton's equations \eqref{Discrete Right Eq} & The discrete right Hamilton's equations \eqref{Discrete Right Eq}\\
		yield the explicit variational integrator & yield the explicit variational integrator \\
		$\begin{aligned}
			r_1 & = r_0 - hCp(q_0^t)^{2p-1} \nabla f(q_0) ,\\
			r_1^t & = r_0^t + h \frac{p(p+1)}{2(q_0^t)^{p+2}} r_1^\top r_1 \\ & \quad \qquad - hCp(2p-1)(q_0^t)^{2p-2} f(q_0),\\
			q_1 &= q_0 +h \frac{p}{(q_0^t)^{p+1}} r_1 , \\
			q_1^t  & = q_0^t + h.
		\end{aligned}$
		& 
		$\begin{aligned}
			r_1 & = r_0 -  \frac{p^2}{\mathring{p}} hC (q_0^t)^{2p-\frac{\mathring{p}}{p}} \nabla f(q_0), \\
			r_{1/2}^t & =  \frac{p^3+ p \mathring{p} }{2 \mathring{p} (q_0^t)^{p+\frac{\mathring{p}}{p} +1} } h  r_1^\top  r_1 + \frac{p \mathring{p} - 2p^3 }{ \mathring{p}  (q_0^t)^{\frac{\mathring{p}}{p} +1-2p}  } hCf(q_0), \\
			r_1^t & = \left[ 1- h (q_0^t)^{- \frac{\mathring{p}}{p}} \left(1-\frac{p}{\mathring{p}}\right)  \right]^{-1}\left(  r_0^t +  r_{1/2}^t  \right)        ,\\
			q_1 &= q_0 +  \frac{p^2}{\mathring{p} } h(q_0^t)^{-p-\frac{\mathring{p}}{p}} r_1,  \\
			q_1^t  & = q_0^t + \frac{p}{\mathring{p} } h(q_0^t)^{1-\frac{\mathring{p}}{p}}.
		\end{aligned}$
	\end{tabular}}
\end{table}

\hfill 

%\linenumbers

Other types of first-order variational integrators can be constructed for the Poincar\'e transformed Hamiltonians, such as Prolongation-Collocation variational integrators \cite{LeSh2011}, Galerkin variational integrators \cite{LeZh2011}, and higher order HTVIs.  We will not consider these integrators here since they require that one solves systems of nonlinear equations and cannot be implemented explicitly in general.  Having said that, in practice, implicit methods for the numerical solution of ODEs that can be solved using fixed-point iterations (as opposed to Newton iterations) can be quite competitive, as there is a good initial guess which may allow them to converge in a small number of iterations, and the iterations are inexpensive as they do not require the assembly of a Jacobian. The convergence of the fixed-point iteration depends on the conditioning of the system of equations, and may impose a stringent time-step restriction. This can be overcome by the use of exponential integrators~\cite{HoOs2010}, and in particular, symplectic and energy-preserving exponential integrators~\cite{ShLe2019}.  

We have also implemented non-variational methods based on these Direct and Adaptive Approaches (the Classical 4th Order Explicit Runge--Kutta Method, and MATLAB's explicit adaptive ODE solvers (\texttt{ode23}, \texttt{ode45})), and more traditional optimization methods have been tested as well such as Nesterov's Accelerated Gradient \eqref{NesterovUpdate} and adaptive optimization algorithms. \\

\noindent \underline{\textbf{Non-Variational Symplectic Integrators based on the Direct \& Adaptive Approaches}}:\\

\begin{enumerate}[label=(\roman*)]
	
	\item \textbf{Direct and Adaptive Approaches with Splitting of the Hamiltonian}: \\
	The Direct Approach with splitting of the Hamiltonian is the approach presented in \cite{JordanSymplecticOptimization}. The three terms of the Poincar\'e transformed Hamiltonian (\ref{pHamiltonian}) are considered separately. They generate dynamics in the extended phase-space via six vector fields, and a symmetric leapfrog composition of the corresponding component dynamics is constructed to obtain a symplectic integrator (referred to in the numerical results section as ``Direct Splitting''). A new symplectic integrator can also be obtained by adapting the approach presented in \cite{JordanSymplecticOptimization} to the adaptive Poincar\'e transformed Hamiltonian (\ref{ADHamiltonian}) to obtain a symplectic integrator (referred to in the numerical results section as ``Adaptive Splitting''). \\
	
	\begin{table}[!h]
		\begin{tabular}{l|l}
			\textbf{Direct Approach} & \textbf{Adaptive Approach}  \\ \hline
			\small
			$
			\begin{aligned}
				t &= t + \frac{h}{2}, \\
				r^t &= r^t + \frac{h}{2}  \frac{p(p+1)}{2t^{p+2}} r^\top r   \\ & \qquad  - \frac{h}{2}  Cp(2p-1) t^{2p-2} f(q), \qquad \\
				r &= r - \frac{h}{2} C p t^{2p-1} \nabla f(q),  \\
				q  &= q + h \frac{p}{t^{p+1}} r, \\
				r  &= r - \frac{h}{2} C p t^{2p-1} \nabla f(q), \\
				r^t &= r^t + \frac{h}{2}  \frac{p(p+1)}{2t^{p+2}} r^\top r   \\ & \qquad  - \frac{h}{2}  Cp(2p-1) t^{2p-2} f(q), \\
				t  &= t + \frac{h}{2}.
			\end{aligned} $
			&           
			\small
			$
			\begin{aligned}
				\quad t &= \left(  t^{\mathring{p}/p} + \frac{h}{2}  \right)^{p/\mathring{p}}, \\
				\theta &= \frac{\mathring{p}}{\mathring{p}-p}  \left[  	\left(  \frac{p^3}{2\mathring{p}} + \frac{p}{2} \right) t^{-p-1} r^\top r + \left( p - \frac{2p^3}{\mathring{p}} \right)  t^{2p-1} C f(q) \right], \\
				r^t  &= (r^t + \theta) \exp{\left(  \left( 1-\frac{p}{\mathring{p}} \right) \frac{h}{2} t^{-\mathring{p}/p} \right)} - \theta ,\\
				r &= r - h\frac{Cp^2}{2 \mathring{p}} t^{2p-\mathring{p}/p} \nabla f(q) , \\
				q  &= q + h \frac{p^2}{\mathring{p} t^{p+\mathring{p}/p}} r ,\\
				r  &=  r - h\frac{Cp^2}{2\mathring{p}} t^{2p-\mathring{p}/p} \nabla f(q), \\
				\theta &= \frac{\mathring{p}}{\mathring{p}-p}  \left[  	\left(  \frac{p^3}{2\mathring{p}} + \frac{p}{2} \right) t^{-p-1} r^\top r + \left( p - \frac{2p^3}{\mathring{p}} \right)  t^{2p-1} C f(q) \right], \\
				r^t  &= (r^t + \theta) \exp{\left(  \left( 1-\frac{p}{\mathring{p}} \right) \frac{h}{2} t^{-\mathring{p}/p} \right)} - \theta,\\
				t  &=  \left(  t^{\mathring{p}/p} + \frac{h}{2}  \right)^{p/\mathring{p}}.
			\end{aligned} $
		\end{tabular}
		\caption{Updates for Direct and Adaptive Approaches with Splitting of the Hamiltonian}
	\end{table}

	\item \textbf{Phase-space Cloning and Splitting}: \\
	\noindent A very natural approach to integrate these non-separable Hamiltonian dynamics consists in defining a new Hamiltonian via two copies of the Poincar\'e transformed Hamiltonian in an extended phase space of dimension twice as large \cite{Piha2015}: 
	\begin{equation}
		\tilde{H}(\bar{q},\tilde{\bar{q}} , \bar{r}, \tilde{\bar{r}} ) = \bar{H}_1(\bar{q} , \tilde{\bar{r}}) + \bar{H}_2(\tilde{\bar{q}},\bar{r} ) ,
	\end{equation}
	where $\bar{H}_1 = \bar{H}_2 = \bar{H}$. Hamilton's equations are then given by
	\begin{equation}
		\dot{\bar{q}} = \nabla_{\bar{r}} \bar{H}_2,   \qquad   \dot{\tilde{\bar{q}}} = \nabla_{\tilde{\bar{r}}} \bar{H}_1,     \qquad  \dot{\bar{r}} = - \nabla_{\bar{q}} \bar{H}_1,   \qquad   \dot{\tilde{\bar{r}}} = \nabla_{\tilde{\bar{q}}} \bar{H}_2. 
	\end{equation}
	We can then integrate this new Hamiltonian system explicitly using a Strang Splitting or a Yoshida 4 or Yoshida 6 Splitting for instance (referred to as ``CloningStrang'', ``CloningY4'', and ``CloningY6'' in the numerical results section). This approach will usually require a larger number of evaluations of the objective function $f$ and of its gradient at each step. 
	
\end{enumerate}

\hfill

\subsection{Numerical Results}

The numerical methods presented in the previous section have been conducted to minimize the quartic function
\begin{equation}\label{Quartic}
	f(x) = \left[(x-1)^\top  \Sigma (x-1) \right]^2,
\end{equation} 
where $x\in \mathbb{R}^{50}$ and  $\Sigma_{ij} =0.9^{|i-j|}$. This convex function achieves its minimum value $0$ at $x^*=1$. \\

Unless specified otherwise, the termination criterion used was 
\begin{equation}\label{TerminationCriterion}
	|f(x_k) - f(x_{k-1})|<\delta  \quad \text{and} \quad \Vert  \nabla f(x_k) \Vert < \delta  \qquad \text{where} \quad \delta = 10^{-10}. 
\end{equation}

\hfill 

\subsubsection{Adaptive versus Direct approach}  

\begin{figure}[!h]
	\hspace*{-25mm}
	\centering
	\begin{minipage}[b]{0.56\textwidth}
		\includegraphics[width=\textwidth]{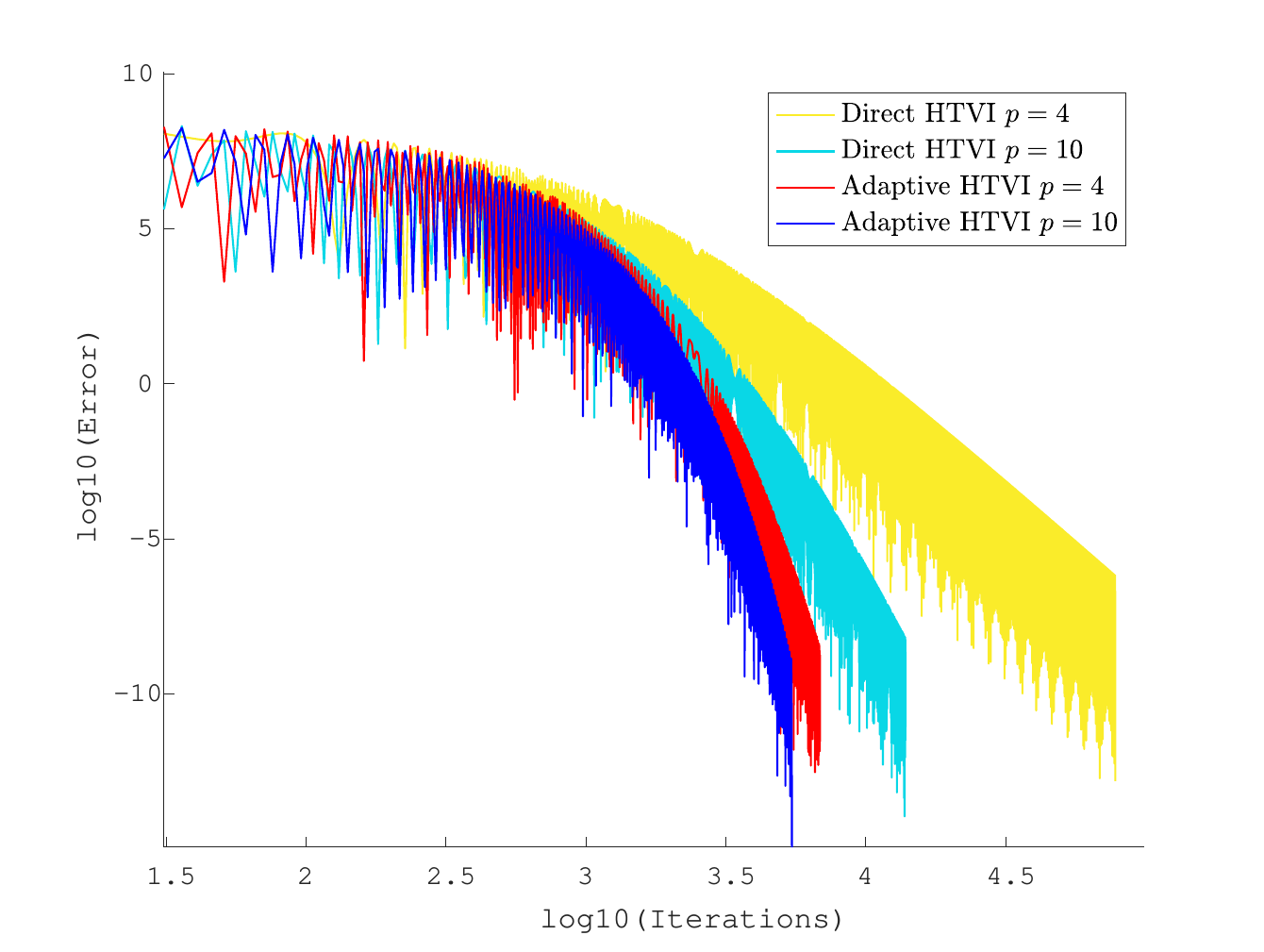}
	\end{minipage}
	\hspace{-10mm}
	\begin{minipage}[b]{0.56\textwidth}
		\includegraphics[width=\textwidth]{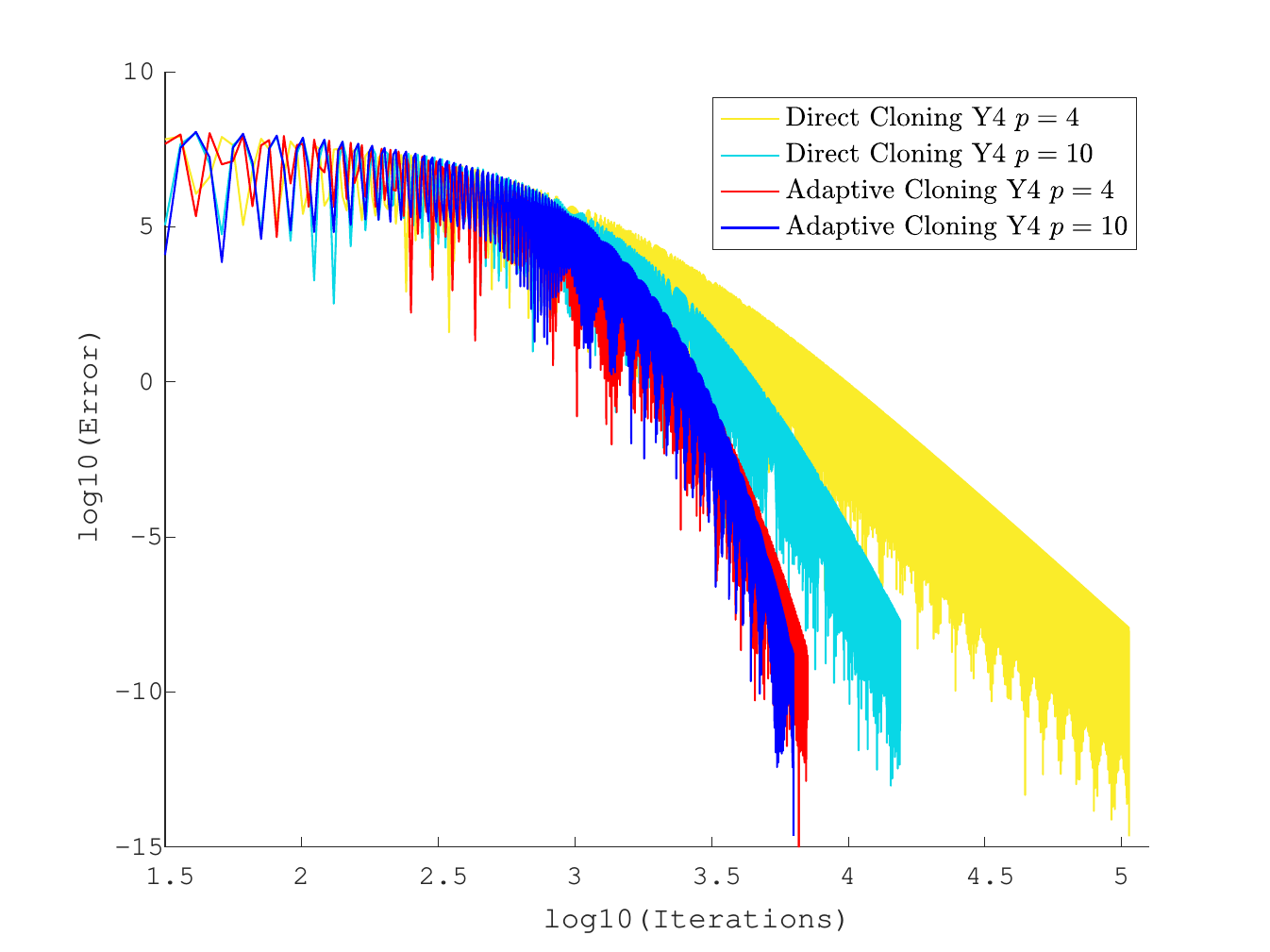}
	\end{minipage}
	\hspace*{-28mm}
	\caption{Comparison of the rates of convergence between the Direct and Adaptive approach for the HTVI method and the Cloning method with a Yoshida 4 splitting. We can clearly see that the Adaptive approach outperforms the Direct approach. \label{Fig:ADvsNonAD}}
\end{figure}

Numerical experiments conducted with all the symplectic algorithms presented in Section \ref{SectionMethods} showed that a carefully tuned Adaptive approach enjoys a significantly better rate of convergence and a much smaller number of steps required to achieve convergence than the Direct approach, as can be seen in Figure \ref{Fig:ADvsNonAD} and Table \ref{Table:ADvsNonAD} for the HTVI and CloningY4 methods. Although the Adaptive approach requires a smaller fictive time-step $h$ than the Direct approach, the physical time steps resulting from $t = \tau^{p/\mathring{p}}$ in the Adaptive approach grow rapidly to values larger than the physical time-step of the Direct Approach.

\begin{table}[!h]
	\resizebox{15cm}{!} {
	\begin{tabular}{lcccc|lcccc}
		\multicolumn{1}{l}{\textbf{Approach}} & \textbf{$p$} & \multicolumn{1}{c}{\textbf{$\mathring{p}$}} & \multicolumn{1}{c}{\textbf{$h$}} & \multicolumn{1}{c|}{\textbf{Iterations}} & \multicolumn{1}{l}{\textbf{Approach}} & \textbf{$p$} & \multicolumn{1}{c}{\textbf{$\mathring{p}$}} & \multicolumn{1}{c}{\textbf{$h$}} & \multicolumn{1}{c}{\textbf{Iterations}} \\ \hline
		Direct HTVI                           & 4            & -                                           & 8.00E-04                         & 77 878                                   & Direct CloneY4                     & 4            & -                                           & 9.00E-04                         & 106 530                                 \\
		Adap. HTVI                         & 4            & 0.5                                         & 1.21E-04                         & 6 867                                    & Adap. CloneY4                   & 4            & 0.5                                         & 1.15E-04                         & 7 101                                   \\
		Direct HTVI                           & 10           & -                                           & 4.00E-04                         & 13 872                                   & Direct CloneY4                     & 10           & -                                           & 3.30E-04                         & 15 498                                  \\
		Adap. HTVI                         & 10           & 0.5                                         & 1.95E-05                         & 5 564                                    & Adap. CloneY4                   & 10           & 0.5                                         & 1.62E-05                         & 6 300                                   \\ \hline
	\end{tabular}}
	\caption{Comparison of the Direct and Adaptive approach for the HTVI method and the Cloning method with a Yoshida 4 splitting. The Adaptive approach clearly outperforms the Direct approach in terms of number of iterations required.   \label{Table:ADvsNonAD} }
\end{table}

\subsubsection{Comparison of Methods within the Direct and Adaptive approaches} 

\begin{figure}[!h]
	\hspace*{-25mm}
	\centering
	\begin{minipage}[b]{0.55\textwidth}
		\includegraphics[width=\textwidth]{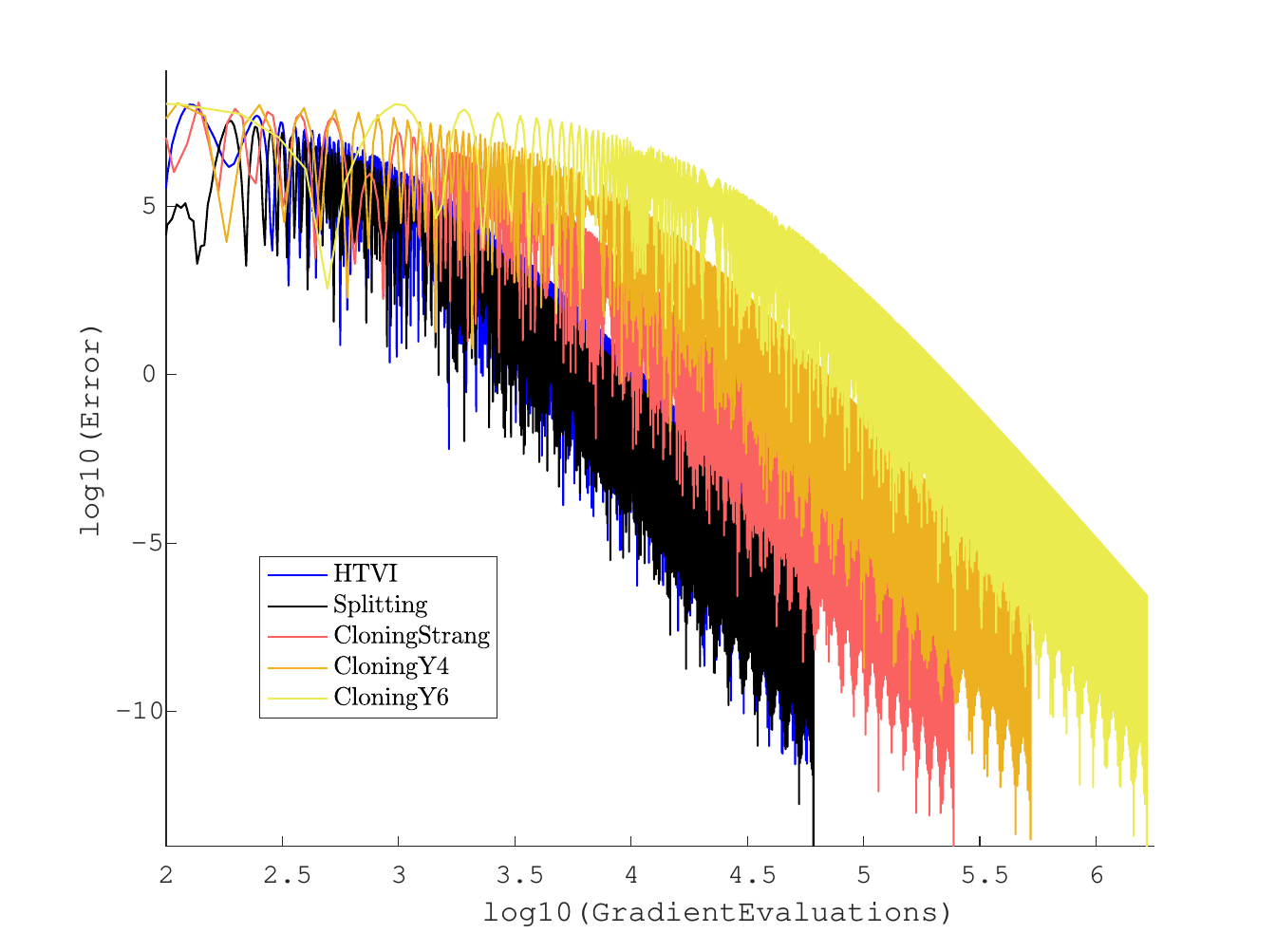}
	\end{minipage}
	\hspace{-10mm}
	\begin{minipage}[b]{0.55\textwidth}
		\includegraphics[width=\textwidth]{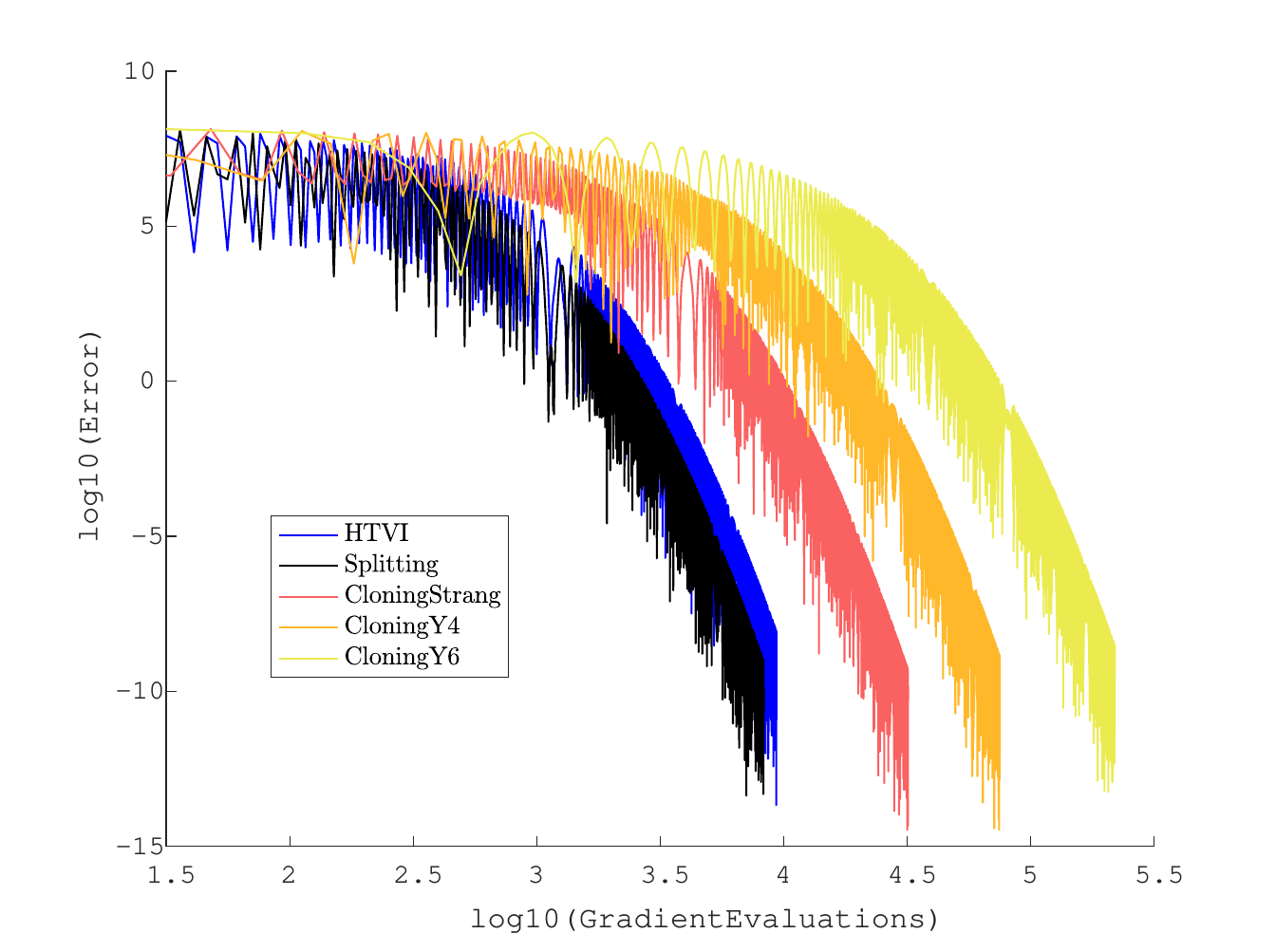}
	\end{minipage}
	\hspace*{-28mm}
	\caption{Comparison of the convergence, in terms of gradient evaluations needed, of the different symplectic integrators within the Direct approach (on the left), and within the Adaptive approach (on the right).  \label{Fig:MethodCompare}}
\end{figure}

Numerical experiments were conducted to compare the various algorithms presented in Section \ref{SectionMethods}, and the results are presented in Figure \ref{Fig:MethodCompare} and Table \ref{Table:MethodCompare}. Although the number of iterations for all methods were of the same order of magnitude, the HTVI method and the Splitting method based on the idea of  \cite{JordanSymplecticOptimization} performed much better than the methods based on the phase-space Cloning idea of \cite{Piha2015}. This is mostly due to the fact that these phase-space Cloning methods require several evaluations of the objective function $f$ and of its gradient $\nabla f$ at each iteration  (3 for Strang Splitting, 7 for Yoshida's 4th Order Splitting, and 19 for Yoshida's 6th Order Splitting), while the HTVI and Splitting methods only required one such evaluation at each iteration. As a result, these phase-space Cloning methods also required much more computational time to achieve convergence. It might be possible to improve the performance of these phase-space Cloning methods by adding an extra term in the final Hamiltonian which binds the two copies of the Poincar\'e transformed Hamiltonian, as was done in \cite{Tao2016}. However, this additional term is likely to require a larger number of compositions when splitting the final Hamiltonian, which would require more gradient evaluations of the objective function at each step. Thus, even though the trick presented in \cite{Tao2016} could reduce the number of iterations required to achieve convergence, it seems very unlikely that the resulting algorithm would be competitive against the HTVI and the Splitting methods, in terms of computational time and total number of gradient evaluations needed.

\begin{table}[!h]
		\resizebox{15cm}{!} {
	\begin{tabular}{lccc|lcccc}
		\textbf{Method}    & \textbf{$p$} & \textbf{$h$} & \textbf{Iterations} & \textbf{Method}      & \textbf{$p$} & \textbf{$\mathring{p}$} & \textbf{$h$} & \textbf{Iterations} \\ \hline
		Direct HTVI        & 4            & 8.7E-04      & 57 504              & Adaptive HTVI        & 4            & 1                       & 2.4E-04      & 9 361               \\
		Direct Splitting   & 4            & 9.5E-04      & 60 881              & Adaptive Splitting   & 4            & 1                       & 2.9E-04      & 8 313               \\
		Direct CloneStrang & 4            & 9.7E-04      & 81 367              & Adaptive CloneStrang & 4            & 1                       & 2.4E-04      & 10 638              \\
		Direct CloneY4     & 4            & 8.9E-04      & 74 747              & Adaptive CloneY4     & 4            & 1                       & 2.9E-04      & 10 721              \\
		Direct CloneY6     & 4            & 7.9E-04      & 87 075              & Adaptive CloneY6     & 4            & 1                       & 2.0E-04      & 11 549              \\ \hline
	\end{tabular}}
	\caption{Number of iterations needed until convergence of the different symplectic integrators within the Direct (on the left) and Adaptive (on the right) approaches.  \label{Table:MethodCompare}}
\end{table}

\subsubsection{Dependence on $p$ and $\mathring{p}$ in the Adaptive approach} 
\begin{figure}[!h]
	\hspace*{-25mm}
	\centering
	\begin{minipage}[b]{0.55\textwidth}
		\includegraphics[width=\textwidth]{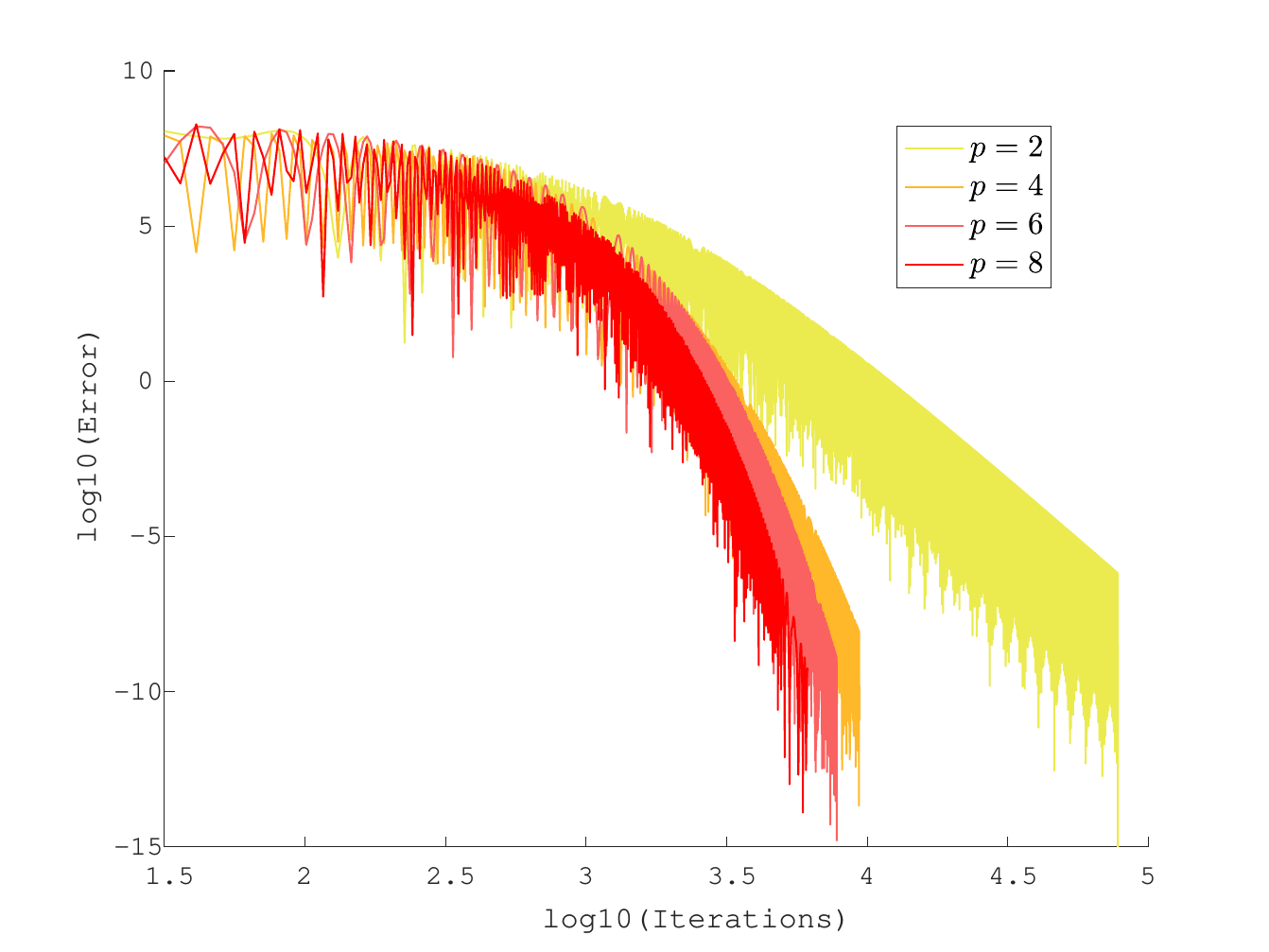}
	\end{minipage}
	\hspace{-10mm}
	\begin{minipage}[b]{0.55\textwidth}
		\includegraphics[width=\textwidth]{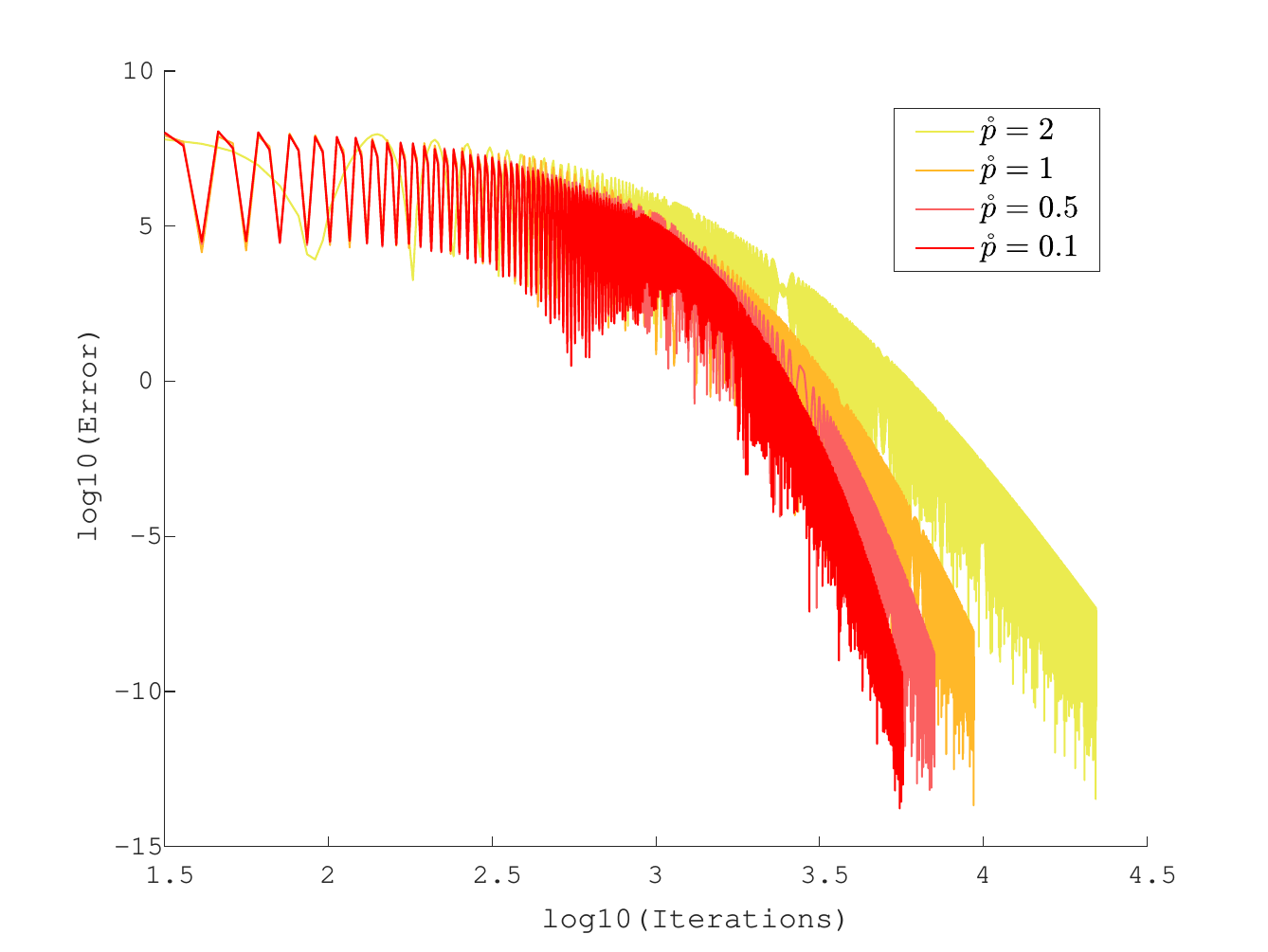}
	\end{minipage}
	\hspace*{-25mm}
	\caption{Evolution of the rates of convergence of the HTVI method as $p$ is increased (on the left), and as $\mathring{p}$ is decreased (on the right).  \label{Fig:Evolution}}
\end{figure}

We conducted numerical experiments with the HTVI method to study the evolution of the performance of the Adaptive approach as the parameters $p$ and $\mathring{p}$ are varied. We can see from the results presented in Figure \ref{Fig:Evolution} and Table \ref{Table:Evolution} that the Adaptive HTVI method becomes more and more efficient as $p$ is increased and $\mathring{p}$ is decreased. The improvement in efficiency is very important as we increase $p$ from $p=2$ to $p=4$, while it is minor but still noticeable as we increase $p$ from $p=4$ to $p=8$. A possible explanation for this behavior is that the integrator might not be of high enough order to distinguish between the $p=6$ and $p=8$ Bregman dynamics. Note that the fictive time-step $h$ must be reduced as $p$ increases or $\mathring{p}$ decreases, but the time relation $t = \tau^{p/\mathring{p}}$ ensures that the resulting physical time steps do not become significantly smaller.

\begin{table}[!h]
	\resizebox{15cm}{!} {
	\begin{tabular}{lcccc|lcccc}
		\textbf{Method} & \textbf{$p$} & \textbf{$\mathring{p}$} & \textbf{$h$} & \textbf{Iterations} & \textbf{Method} & \textbf{$p$} & \textbf{$\mathring{p}$} & \textbf{$h$} & \textbf{Iterations} \\ \hline
		Adaptive HTVI   & 2            & 1                       & 8.0E-04      & 77 855              & Adaptive HTVI   & 4            & 2                       & 3.8E-04      & 22 128              \\
		Adaptive HTVI   & 4            & 1                       & 2.4E-04      & 9 361               & Adaptive HTVI   & 4            & 1                       & 2.4E-04      & 9 361               \\
		Adaptive HTVI   & 6            & 1                       & 9.4E-05     & 7 785               & Adaptive HTVI   & 4            & 0.5                     & 1.2E-04      & 7 099               \\
		Adaptive HTVI   & 8            & 1                       & 6.1E-05     & 6 133               & Adaptive HTVI   & 4            & 0.1                     & 2.4E-05     & 5 689               \\ \hline
	\end{tabular}}
	\caption{Evolution of the fictive time-step $h$ and number of iterations until convergence for the HTVI method as $p$ increases (left), and as $\mathring{p}$ decreases (right).   \label{Table:Evolution} }
\end{table}

\subsubsection{Comparison to non-symplectic integrators}
We will now present the results of numerical experiments investigating the role that symplecticity plays when integrating the Bregman dynamics in the Direct and Adaptive approaches.

We first implemented fixed time-step integrators such as the 4th-order explicit Runge--Kutta method, but these failed to converge both in the Direct and Adaptive approaches.  The reason why convergence cannot be achieved may have to do with the nonautonomous aspect of the differential equation. More precisely, explicit Runge--Kutta methods are conditionally stable, where stability intervals for the time steps depend on the expansivity of the differential equation. Since the differential equations considered here are not autonomous, the stability intervals are time-dependent, and thus any fixed choice of time-step may eventually violate the stability condition. It might be possible to achieve low accuracy convergence using these methods, but the fact that they cannot achieve higher accuracy and are likely to lose stability eventually makes them undesirable.

\begin{figure}[!h]
	\hspace*{-25mm}
	\centering
	\begin{minipage}[b]{0.55\textwidth}
		\includegraphics[width=\textwidth]{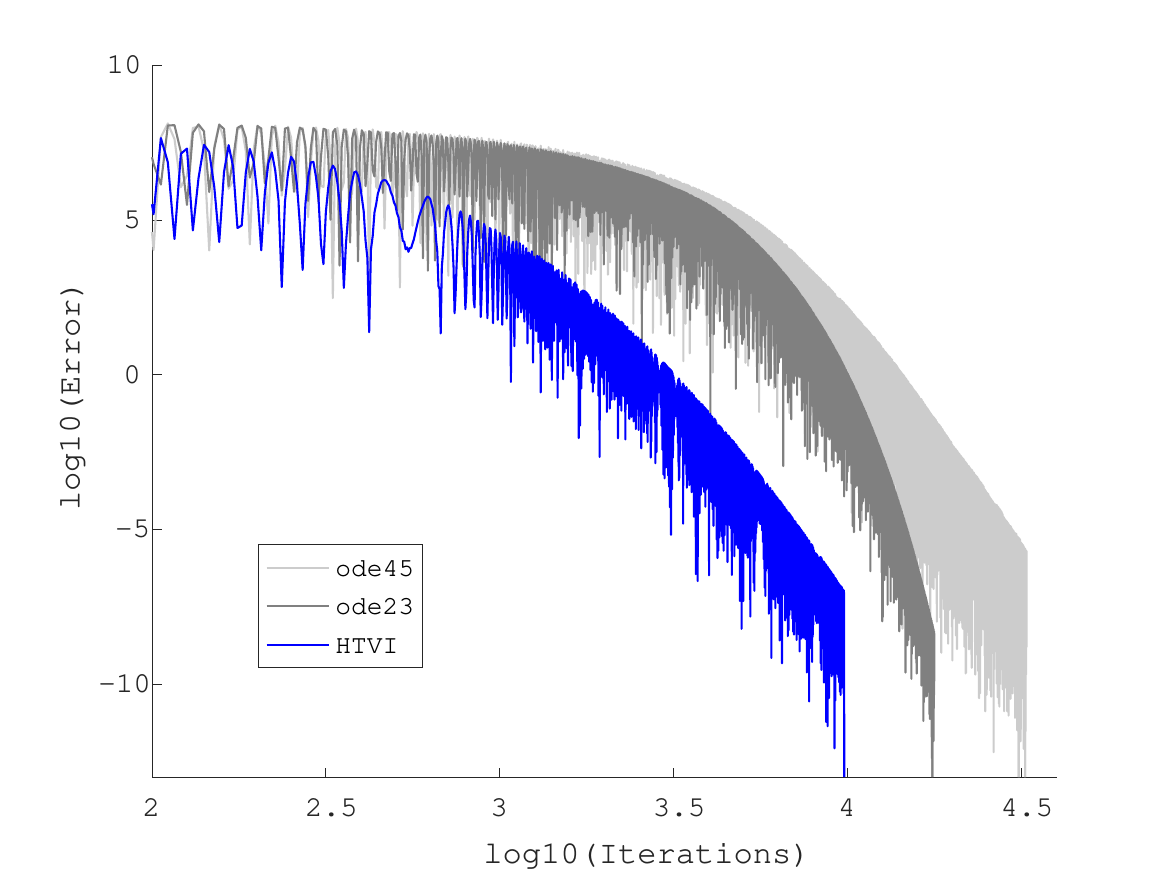}
	\end{minipage}
	\hspace{-10mm}
	\begin{minipage}[b]{0.55\textwidth}
		\includegraphics[width=\textwidth]{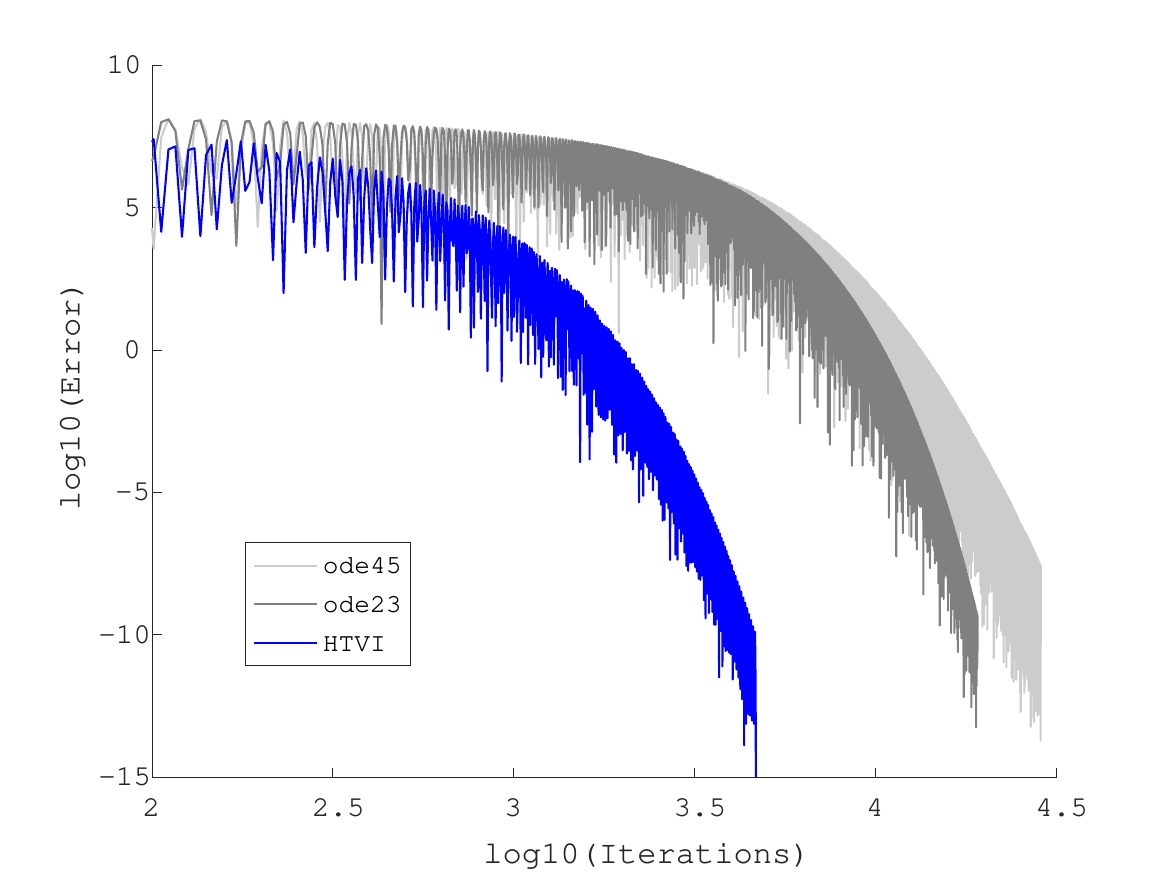}
	\end{minipage}
	\hspace*{-28mm}
	\caption{ Comparison of the HTVI method with the \texttt{ode23} and \texttt{ode45} MATLAB functions in the Direct (left) and Adaptive (right) approaches with $p=10$ and $\mathring{p} = 0.5$. The HTVI method requires outperforms the MATLAB solvers.   \label{Fig:SympTest}}
\end{figure}

We then considered variable time-step explicit Runge--Kutta methods. To this end, we tested MATLAB's differential equation solvers \texttt{ode45} and \texttt{ode23}, which are explicit variable time-step Runge--Kutta pairs, and the corresponding numerical results are presented in Figure \ref{Fig:SympTest}. The HTVI method required a significantly smaller number of iterations than the MATLAB solvers. Furthermore, an inherent part of the time-step control in embedded Runge--Kutta methods is that, at each iteration, the underlying Runge--Kutta method may be executed several times to determine the appropriate time-step that satisfies the prescribed tolerances. For this reason, the MATLAB solvers require more evaluations of $f$ and $\nabla f$ at each iteration, and since they also required more iterations than the HTVI method, these MATLAB solvers are much less competitive. 

It should also be noted that the MATLAB solvers did not exhibit any improvements when used with the Adaptive approach instead of the Direct approach, while the HTVI method improved significantly. This is not surprising since the MATLAB solvers \texttt{ode23} and \texttt{ode45} both use a variable time-step strategy, regardless of the approach chosen.

Note that our Adaptive approach and the embedded Runge--Kutta methods use adaptivity in two fundamentally different ways. Our approach uses \textit{a priori} adaptivity based on known global properties of the family of differential equations considered (i.e. the time-translation symmetry of the family of Bregman dynamics). In contrast, embedded Runge-Kutta methods use adaptivity based on \textit{a posteriori} local error estimates. This could explain why the embedded Runge--Kutta methods do not perform as well as our Adaptive approach: \textit{a posteriori} estimators might focus mostly on the fast local oscillations of the Bregman dynamics and not on the slower global decay, and these fast oscillations might be forcing the embedded Runge--Kutta methods to adaptively take smaller time-steps than necessary. 

We can also see from Figure \ref{Fig:SympTest} that for both the symplectic and non-symplectic adaptive methods, a significant number of iterations are needed before the error effectively starts decaying. The fact that this slow initial behavior persists with those two approaches, which use time-adaptivity in the two fundamentally different ways described in the previous paragraph, suggests that this behavior might be intrinsic to the continuous trajectory being discretized and that time-adaptivity might not be able to help accelerate this initial phase.

\subsubsection{Comparison to popular optimization methods}

Finally, we have compared the performance of our Adaptive HTVI method to Nesterov's Accelerated Gradient (NAG)  \eqref{NesterovUpdate} with the same initial time-step $h=2\times 10^{-6}$, and to popular adaptive optimization algorithms such as Trust Region Steepest Descent (TRUST), ADAM \cite{ADAM}, AdaGrad \cite{AdaGrad}, and RMSprop \cite{RMSprop}. 

\begin{figure}[!h]
	\includegraphics[width=1\textwidth]{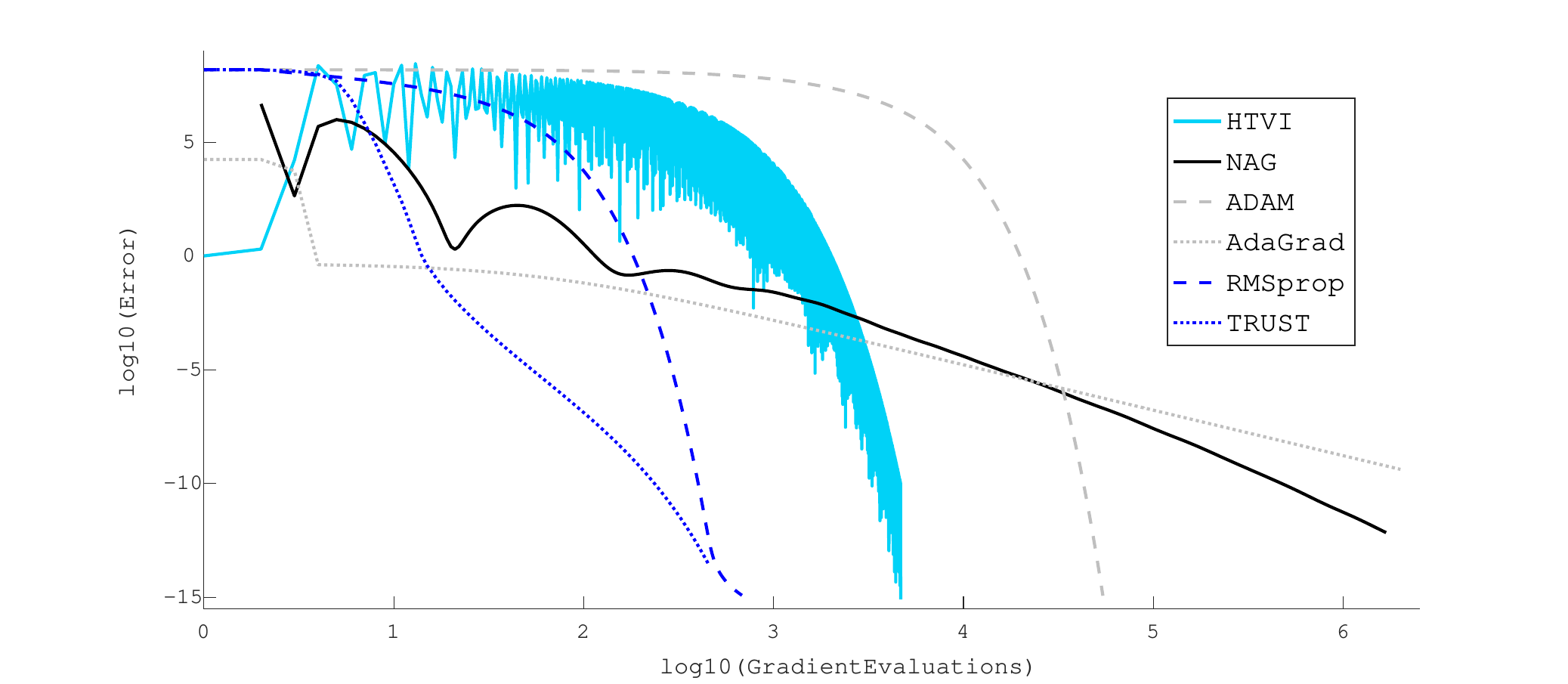}
	\caption{Comparison of HTVI, NAG, and other adaptive optimization algorithms to achieve convergence on the quartic objective function \eqref{Quartic}, with different values of $\delta$ as termination criterion \eqref{TerminationCriterion}. Note that HTVI and NAG were implemented with the same initial time step.  \label{Fig:OptMethods}}
\end{figure}

Figure \ref{Fig:OptMethods} and Table \ref{table3} present the numerical results obtained when applying these algorithms to the quartic objective function \eqref{Quartic}. Although the Adaptive HTVI method is not the most efficient method, we can see that it significantly outperformed certain popular optimization algorithms on this particular convex problem. This suggests that the Adaptive HTVI method might be a competitive first-order explicit algorithm, and that it might be worth considering it as one of several possible options to use in practice, as the relative performance often depends on the specific choice of objective function.

\begin{table}[ht]
	\centering
	\resizebox{0.95\textwidth}{!}{\begin{tabular}{r|cccccccccc}
			\hline
			$\delta = $ & $ 10^{-2} $ &  $ 10^{-3} $ & $10^{-4} $ & $10^{-5} $ & $ 10^{-6}  $ & $ 10^{-7}  $ & $ 10^{-8}  $ & $ 10^{-9} $ & $ 10^{-10}$ \\
			\hline
			HTVI & 2 182  &  2 486 &	2 750 &	3 233  &	3 434 & 3 593	&  4 014	& 4 097 &	4 566\\
			NAG &  4 143  & 10 949  & 27 660	 & 65 724	& 154 258   & 341 928  & 745 292  & 1.7E6	& $>$1E10 \\
			ADAM &  29 665  & 32 733  & 35 802	 & 38 871	& 41 939   & 45 008  & 48 076 & 51 145	& 54 215 \\
			AdaGrad &  520 482  & 2.4E06  & 1.1E07 	 &  5.2E07	&  2.4E08   & --- & --- & ---	& --- \\
			RMSprop &  276  & 305    &  334	 & 	363 &  393   &  422 &  452 &  498	& 682  \\
			TRUST &  32  & 48  & 71	 & 106	& 154   & 215  & 288 & 366	& 455 \\
			\hline
	\end{tabular}}
	\caption{Comparison of the number of iterations needed for HTVI, for NAG, and for other adaptive optimization algorithms to achieve convergence on the quartic objective function \eqref{Quartic}, with different values of $\delta$ as termination criterion \eqref{TerminationCriterion}. Note that HTVI and NAG were implemented with the same initial time step. For all these algorithms, the number of gradient evaluations equals the number of iterations.   }	\label{table3}
\end{table}

\begin{remark*}
	In \cite{JordanSymplecticOptimization}, the authors noted that Nesterov's Accelerated Gradient algorithm transitions into an exponential rate of convergence once it is sufficiently close to the minimum of certain objective functions, and suggested that this behavior requires strong convexity of the objective function in the neighborhood of the minimum. Similarly to the strategy presented in \cite{JordanSymplecticOptimization}, a gradient flow can be incorporated into the updates of the Direct and Adaptive algorithms presented so that for certain objective functions, the same exponential rate of convergence can be achieved close to the minimum.
\end{remark*}

\begin{remark*}
	In very high-dimensional nonconvex optimization problems of practical interest, it has been noted empirically that a main source of difficulty is not the presence of local minima but rather the ubiquity of saddle points surrounded by high error plateaux \cite{Dauphin2014,Jin2021}. These saddle points can significantly slow down gradient-based algorithms and give the illusion of the existence of a local minima. It was demonstrated in \cite{Jin2018} via a variant of Nesterov's Accelerated Gradient algorithm that momentum techniques can escape saddle points faster than standard gradient methods and can thereby accelerate convergence in the nonconvex setting as well. This suggests that the variational framework for accelerated optimization and our Adaptive approach to obtain symplectic optimization algorithms may also be promising with regards to nonconvex optimization. 
\end{remark*}

\section{Conclusions}

Due to the degeneracy of the Hamiltonian, adaptive variational integrators based on the Poincar\'e transformation should be constructed using discrete Hamiltonians, which are Type II or III generating functions. This has potential implications for the numerical properties of such integrators, and might explain why there has only been a limited amount of work on the construction of adaptive variational integrators based on the traditional Lagrangian perspective. The efficiency of the resulting integrator is largely based upon a proper choice of the monitor function $g$, and more research is needed to find a general choice of $g$ that maintains a decent level of efficiency.

We have also noted that the gain in efficiency provided by adaptivity depends on the properties of the Hamiltonian dynamical system, and tends to be more significant in regions of high curvature of the Hamiltonian in the vector field. We focused primarily on Taylor variational integrators, but Galerkin variational integrators are likely to be very promising as well, since the cost of evaluating the monitor function and its derivatives should be low. In addition, the Galerkin approximation scheme may help inform a better choice of monitor function, due to the extensive literature on efficient \textit{a posteriori} error estimation. \textit{A posteriori} error estimation, in general, would be a nice addition to give some guarantees on global accuracy.

Finally, we used our adaptive framework together with the variational approach to accelerated optimization presented in \cite{WiWiJo16} to design efficient Hamiltonian variational and non-variational explicit integrators for symplectic accelerated optimization. We noted that a careful use of adaptivity and symplecticity can result in significantly faster algorithms, and it could also play an important role for nonconvex optimization problems. It would be desirable to understand at a theoretical level the role that adaptivity and symplecticity plays in the accurate and stable discretization of flows that correspond to accelerated optimization algorithms, which could better inform the choice of monitor functions, and the convex function used to define the Bregman divergence that arises in the construction of the Bregman Hamiltonian. This direction seems particularly promising for constructing novel optimization algorithms with superior computational efficiency and performance. Another possible future research direction is to consider how these variational and adaptive frameworks extend to more general spaces such as Lie groups \cite{Tao2020,Lee2021} and Riemannian manifolds \cite{Duruisseaux2021Riemannian,Duruisseaux2021Constrained,Duruisseaux2021Projection}. It could also be interesting to consider the implications of this work for stochastic gradient descent methods~\cite{RoMo1951}, by considering it in the context of a Bregman Lagrangian or Hamiltonian, but with a stochastic perturbation of the potential. This naturally leads to considering stochastic generalizations of the adaptive Hamiltonian variational integrators considered in this paper, by extending the existing work on stochastic variational integrators~\cite{BROw2009,HoTy2018}. \\

\section*{Acknowledgements}
The authors were supported in part by NSF under grants DMS-1411792, DMS-1345013, DMS-1813635, by AFOSR under grant FA9550-18-1-0288, and by the DoD under grant 13106725 (Newton Award for Transformative Ideas during the COVID-19 Pandemic).

\appendix
\section{Proof of Theorem \ref{HTVITheorem}}  \label{HTVIProof}

\hfill 

The proof of Theorem \ref{HTVITheorem} is similar to the one presented in the Appendix of  \cite{ ScShLe2017} for Lagrangian Taylor variational integrators. We first start with the right Hamiltonian Taylor variational integrator case. Let $q(t)$ and $p(t)$ denote the solutions of Hamilton's boundary value problem
$$ \dot{q}(t) = g(q(t),p(t),t),  \qquad \dot{p}(t) = f(q(t),p(t),t),  \qquad q(0)=q_0,   \quad p(h)= p_1.$$

Let $q_1 = q(h)$ and $p_0=p(0)$. 
\begin{lemma} \label{Lemma1HTVI}
	Given a  $r$-order Taylor method $\Psi_h^{(r)}$ approximating the exact time-$h$ flow map corresponding to Hamilton's equations, let $\tilde{p}_0$ solve the problem $ p_1 = \pi_{T^*Q} \circ \Psi_h^{(r)}(q_0,\tilde{p}_0).$ Then, 
	$$ \tilde{p}_0 = p_0 + \mathcal{O}(h^{r+1}).$$ 
	\begin{proof}
		Solving the equation $ p_1 = \pi_{T^*Q} \circ \Psi_h^{(r)}(q_0,\tilde{p}_0)  $ for $\tilde{p}_0$  yields 
		\begin{equation*}
			 \tilde{p}_0 = p_1 - \sum_{k=1}^{r}{\frac{h^k}{k!} f^{(k-1)}(q_0,\tilde{p}_0,0)}. 
		\end{equation*} 
		The exact solution $p(t)$ belongs to $C^{r+1}([0,h])$ so Taylor's Theorem gives
		$$ p_0 = p_1 - \sum_{k=1}^{r}{\frac{h^k}{k!} f^{(k-1)}(q_0,p_0,0)} + R_r(h).$$
		Now, since $p(t)$ belongs to $C^{r+1}([0,h])$, $f^{(k-1)}$ is Lipschitz continuous in its arguments for $k=1,...,r-1$. Let $M $ be the largest of the corresponding $(r-1)$ Lipschitz constants with respect to the second argument over the compact interval $[0,C]$. Then, using the triangle inequality,
		\small \begin{equation*}
			\Vert \tilde{p}_0 - p_0 \Vert  = \left\Vert R_r(h) - \sum_{k=1}^{r}{\frac{h^k}{k!}  \left[ f^{(k-1)}(q_0,\tilde{p}_0,0) -  f^{(k-1)}(q_0,p_0,0) \right]  }   \right\Vert   \leq  M \sum_{k=1}^{r}{ \frac{h^k}{k!}  \left\Vert \tilde{p}_0 - p_0 \right\Vert  }  + \Vert R_r(h)  \Vert .
		\end{equation*}
		\normalsize Thus, $ \left(1- M \sum_{k=1}^{r}{ \frac{h^k}{k!}} \right) \Vert \tilde{p}_0 - p_0 \Vert \leq    \Vert R_r(h)  \Vert  = \mathcal{O}(h^{r+1}),$ and by continuity, $\exists \tilde{C} \in (0,C)$ such that $\forall h \in (0,\tilde{C})$, the term $\left(1- M \sum_{k=1}^{r}{ \frac{h^k}{k!}} \right)$ is bounded away from zero, which concludes the proof.
	\end{proof}
\end{lemma}

We now show that starting the $r$-order Taylor method with initial conditions $(q_0,\tilde{p}_0)$ rather than $(q_0,p_0)$ will not affect the order of accuracy of the method.

\begin{lemma}\label{Lemma2HTVI}
	The  $r$-order Taylor method $\Psi_h^{(r)}$ with initial conditions $(q_0,\tilde{p}_0)$ and where $\tilde{p}_0$ solves  $ p_1 = \pi_{T^*Q} \circ \Psi_h^{(r)}(q_0,\tilde{p}_0)  ,$  is accurate to at least $\mathcal{O}(h^{r+1})$ for the Hamiltonian boundary-value problem.
	\begin{proof}
		Let $(\tilde{q}(t),\tilde{p}(t))$ denote the exact solution to Hamiltonian's equations with initial values $(q_0,\tilde{p}_0)$, and let $(q_d(t),p_d(t))$ denote the values generated by the $r$-order Taylor method with initial conditions $(q_0,\tilde{p}_0)$. The Hamiltonian initial-value problem is well-posed, so denoting the Lipschitz constant with respect to the second argument by $M$, we get
				\begin{equation*}
			\begin{aligned}
			\Vert  (q(t),p(t)) - (q_d(t),p_d(t)) \Vert	& \leq  	\Vert  (q(t),p(t)) - (\tilde{q}(t),\tilde{p}(t)) \Vert	 + 	\Vert  (\tilde{q}(t),\tilde{p}(t)) - (q_d(t),p_d(t)) \Vert	\\
			& \leq M \Vert p_0 - \tilde{p}_0  \Vert + \mathcal{O}(h^{r+1})  \leq  \mathcal{O}(h^{r+1}) ,
				\end{aligned}
		\end{equation*}
		where we have used the triangle inequality, and the fact that the local truncation error of a $r$-order Taylor method is $\mathcal{O}(h^{r+1}) $ to bound $\Vert  (\tilde{q}(t),\tilde{p}(t)) - (q_d(t),p_d(t)) \Vert$.
	\end{proof}
\end{lemma}

We are now ready to prove Theorem \ref{HTVITheorem} for right Hamiltonian Taylor variational integrators.

\begin{theorem}
	Consider a Hamiltonian $H$ such that $H$ and $\frac{\partial H}{\partial p}$ are Lipschitz continuous in both variables. Given a $r$-order accurate Taylor method $\Psi_h^{(r)}$ and a $s$-order accurate quadrature formula with weights and nodes $(b_i,c_i)$, define the associated Taylor discrete right Hamiltonian
	$$ H_d^+(q_0,p_1;h) = p_1^\top \tilde{q}_1 - h \sum_{i=1}^{m}{b_i \left[  p_{c_i}^\top  \dot{q}_{c_i} - H(q_{c_i},p_{c_i})  \right]},$$
	where $\tilde{p}_0$ solves $ p_1 = \pi_{T^*Q} \circ \Psi_h^{(r)}(q_0,\tilde{p}_0)$, where $ \tilde{q}_1 = \pi_{Q}  \circ \Psi_h^{(r+1)}(q_0,\tilde{p}_0) $ and $ (q_{c_i},p_{c_i})  = \Psi_{c_i h}^{(r)}(q_0,\tilde{p}_0)$, and where we use the continuous Legendre Transform to obtain $\dot{q}_{c_i}$.
	
	Then, $H_d^+$ approximates $H_d^{+,E}$ with order of accuracy at least $\min{(r+1,s)}$. By Theorem 2.2 in \cite{ScLe2017}, the associated discrete right Hamiltonian map has the same order of accuracy.
	\begin{proof}
		From Lemma \ref{Lemma2HTVI} we have that $ q(c_i h) = q_{c_i} + \mathcal{O}(h^{r+1})$ and $p(c_i h) = p_{c_i} + \mathcal{O}(h^{r+1}), $ and since $\frac{\partial H}{\partial p}$ is Lipschitz in both variables
		$  \dot{q}(c_ih) - \dot{q}_{c_i} =  \frac{\partial H}{\partial p}(q(c_ih),p(c_ih)) -  \frac{\partial H}{\partial p}(q_{c_i},p_{c_i}) =\mathcal{O}(h^{r+1}).  $
		Since the quadrature formula is of order $s$ accurate, equation \eqref{exact_Hd} for $H_d^{+,E} (q_0,p_1 ; h)$ gives
		\begin{equation*}
			H_d^{+,E} (q_0,p_1 ; h)   = p_1^\top  q_1 - h \sum_{i=1}^{m}{b_i \left[  p(c_i h)^\top  \dot{q}(c_ih)  - H\left(q_{c_i} +\mathcal{O}(h^{r+1})  ,p_{c_i} + \mathcal{O}(h^{r+1}) \right)  \right]} + \mathcal{O}(h^{s+1}) .
		\end{equation*}
		Now, since $ \tilde{q}_1 = \pi_{Q}  \circ \Psi_h^{(r+1)}(q_0,\tilde{p}_0) $, it follows from Lemma \ref{Lemma2HTVI} that $ \tilde{q}_1 = q_1 + \mathcal{O}(h^{r+2}) $. Therefore, combining this with the fact that $H$ is Lipschitz continuous in both variables yields
		\begin{equation*}
			\begin{aligned}
			H_d^{+,E} (q_0,p_1 ; h)  & = p_1^\top  \tilde{q}_1 - h \sum_{i=1}^{m}{b_i \left[p_{c_i}^\top  \dot{q}_{c_i}  - H\left(q_{c_i} ,p_{c_i} \right)  \right]} + \mathcal{O}(h^{r+2})  + \mathcal{O}(h^{s+1}) \\ & = H_d^+(q_0,p_1;h) + \mathcal{O}(h^{\min{(r+1,s)} +1}). 
			\end{aligned}
		\end{equation*}
		Therefore, $H_d^+$ approximates $H_d^{+,E}$ with order of accuracy at least $\min{(r+1,s)}$.  
	\end{proof}
\end{theorem}

\hfill 

Theorem \ref{HTVITheorem} can be proven in a similar way for left Hamiltonian Taylor variational integrators. Now, $q(t)$ and $p(t)$ denote the solutions of the Hamilton's boundary-value problem
$$ \dot{q}(t) = g(q(t),p(t),t),  \qquad \dot{p}(t) = f(q(t),p(t),t),  \qquad q(h)=q_1,   \quad p(0)= p_0,$$
and let $q_0 = q(0)$ and $p_1=p(h)$.  Lemma \ref{Lemma1HTVI} is replaced by

\begin{lemma} \label{LeftLemma1HTVI}
	Given a  $(r+1)$-order Taylor method $\Psi_h^{(r+1)}$ approximating the exact time-$h$ flow map corresponding to Hamilton's equations, let $\tilde{q}_0$ solve the problem $ q_1 = \pi_{Q} \circ \Psi_h^{(r+1)}(\tilde{q}_0,p_0).$ Then, 
	$$ \tilde{q}_0 = q_0 + \mathcal{O}(h^{r+2}).$$

	\begin{proof}
		We proceed as in the proof of Lemma \ref{Lemma1HTVI}. We first solve $ q_1 = \pi_{Q} \circ \Psi_h^{(r+1)}(\tilde{q}_0,p_0) $ for $\tilde{q}_0$, and then Taylor expand the exact solution $q(t)$ which belongs to $C^{r+2}([0,h])$. Now, $q(t)$ is Lipschitz continuous in its arguments for $k=1,...,r$, so we can let $M $ be the largest of the corresponding $r$ Lipschitz constants with respect to the first argument over the compact interval $[0,C]$. Then, as before, the triangle inequality can be used to get that
		$ \left(1- M \sum_{k=1}^{r+1}{ \frac{h^k}{k!}} \right) \Vert \tilde{q}_0 - q_0 \Vert  = \mathcal{O}(h^{r+2}) , $ and by continuity, the term inside the parenthesis is bounded away from zero. 
	\end{proof}
\end{lemma} 

In analogy to Lemma \ref{Lemma2HTVI}, we now show that starting the $r$-order Taylor method with initial conditions $(\tilde{q}_0,p_0)$ rather than $(q_0,p_0)$ will not affect the order of accuracy of the method.

\begin{lemma}\label{LeftLemma2HTVI}
	The  $r$-order Taylor method $\Psi_h^{(r)}$ with initial conditions $(\tilde{q}_0,p_0)$ and where $\tilde{q}_0$ solves	$ q_1 = \pi_{Q} \circ \Psi_h^{(r)}(\tilde{q}_0,p_0) ,$  is accurate to at least $\mathcal{O}(h^{r+1})$ for the Hamiltonian boundary-value problem.
	\begin{proof}
		Let $(\tilde{q}(t),\tilde{p}(t))$ denote the exact solution to Hamiltonian's equations with initial values $(\tilde{q}_0,p_0)$, and let $(q_d(t),p_d(t))$ denote the values generated by the $r$-order Taylor method with initial conditions $(\tilde{q}_0,p_0)$. The Hamiltonian initial-value problem is well-posed, so denoting the Lipschitz constant with respect to the first argument by $M$, we get
				\begin{equation*}
		\begin{aligned}
			\Vert  (q(t),p(t)) - (q_d(t),p_d(t)) \Vert	& \leq  	\Vert  (q(t),p(t)) - (\tilde{q}(t),\tilde{p}(t)) \Vert	 + 	\Vert  (\tilde{q}(t),\tilde{p}(t)) - (q_d(t),p_d(t)) \Vert	\\
			& \leq M \Vert q_0 - \tilde{q}_0  \Vert + \mathcal{O}(h^{r+1})  \leq  \mathcal{O}(h^{r+1}) ,  
			\end{aligned}
	\end{equation*}
		where we have used the triangle inequality, and the fact that the local truncation error of a $r$-order Taylor method is $\mathcal{O}(h^{r+1}) $ to bound $\Vert  (\tilde{q}(t),\tilde{p}(t)) - (q_d(t),p_d(t)) \Vert$. 
	\end{proof}
\end{lemma} 

\hfill 

We are now ready to prove Theorem \ref{HTVITheorem} for left Hamiltonian Taylor variational integrators.

\begin{theorem}
	Consider a Hamiltonian $H$ such that $H$ and $\frac{\partial H}{\partial p}$ are Lipschitz continuous in both variables. Given a $r$-order accurate Taylor method $\Psi_h^{(r)}$ and a $s$-order accurate quadrature formula with weights and nodes $(b_i,c_i)$, define the associated Taylor discrete left Hamiltonian
	$$ H_d^-(q_1,p_0;h) = -p_0^\top  \tilde{q}_0 - h \sum_{i=1}^{m}{b_i \left[  p_{c_i}^\top  \dot{q}_{c_i} - H(q_{c_i},p_{c_i})  \right]},$$
	where $\tilde{q}_0$ solve the problem $ q_1 = \pi_{Q} \circ \Psi_h^{(r+1)}(\tilde{q}_0,p_0) , $ where $ (q_{c_i},p_{c_i})  = \Psi_{c_i h}^{(r)}(q_0,\tilde{p}_0),$ and where we use the continuous Legendre Transform to obtain $\dot{q}_{c_i}$. 
	
	Then, $H_d^-$ approximates $H_d^{-,E}$ with order of accuracy at least $\min{(r+1,s)}$.  By a result analogous to Theorem 2.2 in \cite{ScLe2017}, the associated discrete left Hamiltonian map has the same order of accuracy.
	\begin{proof}
		From Lemma \ref{LeftLemma2HTVI} we have that $ q(c_i h) = q_{c_i} + \mathcal{O}(h^{r+1})$, and $p(c_i h) = p_{c_i} + \mathcal{O}(h^{r+1}), $
		and since $\frac{\partial H}{\partial p}$ is Lipschitz in both variables
		$  \dot{q}(c_ih) - \dot{q}_{c_i} =  \frac{\partial H}{\partial p}(q(c_ih),p(c_ih)) -  \frac{\partial H}{\partial p}(q_{c_i},p_{c_i}) =\mathcal{O}(h^{r+1}).  $
		Since the quadrature formula is of order $s$ accurate, equation \eqref{exact_LeftHd} for $H_d^{-,E} (q_1,p_0 ; h)$ gives
		\begin{equation*}
			H_d^{-,E} (q_1,p_0 ; h)  = - p_0^\top q_0  - h \sum_{i=1}^{m}{b_i \left[  p(c_i h)^\top \dot{q}(c_ih)  - H\left(q_{c_i} +\mathcal{O}(h^{r+1})  ,p_{c_i} + \mathcal{O}(h^{r+1}) \right)  \right]} + \mathcal{O}(h^{s+1}) .
		\end{equation*}
		Now, since $ q_1 = \pi_{Q}  \circ \Psi_h^{(r+1)}(\tilde{q}_0,p_0) $, it follows from Lemma \ref{LeftLemma1HTVI} that $ \tilde{q}_0 = q_0 + \mathcal{O}(h^{r+2}) $. Therefore, combining this with the fact that $H$ is Lipschitz continuous in both variables yields
		\begin{equation*}
			H_d^{-,E} (q_1,p_0 ; h)    = H_d^-(q_1,p_0;h) + \mathcal{O}(h^{\min{(r+1,s)} +1}) .
		\end{equation*}
	\vspace{0mm} 
	\end{proof}
\end{theorem}

\newpage

\bibliography{AdaptiveVI}
\bibliographystyle{siamplain}

\end{document}